\documentclass{amsart}
\usepackage{amssymb}

\usepackage{graphicx}
\usepackage{amscd}
\usepackage{scalerel}
\usepackage{stackengine,wasysym}
\usepackage{tikz-cd}
\usepackage{bigints}

\usepackage[utf8]{inputenc}
\usepackage[T1]{fontenc}

\usepackage{enumitem}

\newtheorem{theorem}{Theorem}[section]
\theoremstyle{plain}

\newtheorem{corollary}{Corollary}[section]

\newtheorem{lemma}{Lemma}[section]

\newtheorem{problem}{Problem}
\newtheorem{proposition}{Proposition}[section]
\newtheorem{remark}{Remark}[section]

\numberwithin{equation}{section}

\newcommand{\R}{\mathbb{R}}
\newcommand{\Sp}{\mathrm{Span}}

\begin{document}
\title[Reachable sets]{On the Heisenberg sub-Lorentzian structure on $\R^{2n+1}$ }

\author[M. Grochowski]{Marek Grochowski}
\address{Faculty of Mathematics and Natural Sciences\\ Cardinal Wyszy\'nski University\\
ul. Dewajtis 5, 01-815 Warszawa, Poland}
\email{m.grochowski@uksw.edu.pl}
\date{\today}

\begin{abstract} The aim of this paper is to study the Heisenberg sub-Lorentzian metric on $\R^{2n+1}$, $n\geq 2$. Among other results, we determine its isometry group and provide a basis for the Lie algebra of infinitesimal isometries. We also prove the existence of a maximizing geodesic between any two causally related points and show that every such geodesic is Hamiltonian and hence smooth. Moreover, every maximizing geodesic joining two chronologically related points is timelike. Finally, we describe the boundary of the reachable set and compute the future null and timelike conjugate and cut loci of the origin.
\end{abstract}

\maketitle

\section{Introduction and statement of the results}

Let $M$ be a smooth (i.e., of class $C^{\infty}$) manifold and let $H$ be a smooth bracket-generating distribution of constant rank on $M$. Suppose that $g$ is a smooth pseudo-Riemannian metric on $H$. The pair $(H,g)$ is called a \textit{sub-pseudo-Riemannian structure} on $M$, and the triple $(M,H,g)$ is a \textit{sub-pseudo-Riemannian manifold}. If $g$ is Riemannian, then $(M,H,g)$ is called a \textit{sub-Riemannian manifold}, whereas if $g$ is Lorentzian, then $(M,H,g)$ is called a \textit{sub-Lorentzian manifold}. Sub-Riemannian geometry has been studied extensively for many years, and there is a vast literature devoted to it; see for instance \cite{Agr, AgrIsomp, AgrBar,Bell, leDonne,BoscAgr,Strich, Montg,GroWar,Gro2025} and the reference sections therein. By contrast, sub-Lorentzian geometry is less well developed. Nevertheless, it has attracted increasing interest, and a growing number of papers are devoted to the subject; see for instance \cite{Gro2002,Gro2006,Gro2009,Gro2011,Gro2012,SachkSachk,GodKorMar,Lokut,Grong,Karm}.

In this paper, we study what may be regarded as a model example of a contact sub-Lorentzian manifold. To introduce it, consider the space $\R^{2n+1}$ with coordinates $(x_1,y_1,\dots,x_n,y_n,z)$. Let
\begin{align}\label{fields}
 X_i = \frac{\partial}{\partial x_i} + \frac12 y_i\frac{\partial}{\partial z}, \quad Y_i = \frac{\partial}{\partial y_i} - \frac12 x_i\frac{\partial}{\partial z}\text{,}
\end{align}
$i=1,\dots,n$, and define $H$ to be 
\[
 H = \Sp\{X_1,Y_1,\dots,X_n,Y_n\}\text{.}
\]
It is easy to verify that $H$ is a contact distribution. Indeed, $[X_i,Y_i] = -\frac{\partial}{\partial z}$.
Equip $H$ with a Lorentzian metric $g$ for which the frame $X_1,Y_1,\dots,X_n,Y_n$ is orthonormal, with $X_1$ timelike. We choose $X_1$ as a time orientation of $(\R^{2n+1},H,g)$. 
The pair $(H,g)$ is called the \textit{Heisenberg sub-Lorentzian structure} on $\R^{2n+1}$. 

For the reader’s convenience, we briefly introduce the terminology needed to state our main results. All relevant terminology, facts, and further details are provided in Section \ref{SecPrelimSL}. 
An absolutely continuous curve $\gamma\colon[a,b]\longrightarrow\mathbb{R}^{2n+1}$ is called horizontal if $\dot{\gamma}(t)\in H_{\gamma(t)}$ for almost every $t\in[a,b]$. A horizontal curve $\gamma$ is called timelike future directed if 
\[ 
 g(\dot{\gamma}(t),\dot{\gamma}(t))<0 \quad\text{and}\quad g(\dot{\gamma}(t),X_1(\gamma(t)))<0 
 \] 
 for almost every $t\in[a,b]$. It is called null future directed if 
\[ 
  g(\dot{\gamma}(t),\dot{\gamma}(t))=0, \qquad \dot{\gamma}(t)\ne0, \qquad g(\dot{\gamma}(t),X_1(\gamma(t)))<0 
\] for almost every $t\in[a,b]$. Finally, it is called nonspacelike future directed if 
\[ 
  g(\dot{\gamma}(t),\dot{\gamma}(t))\leq0, \qquad \dot{\gamma}(t)\ne0, \qquad g(\dot{\gamma}(t),X_1(\gamma(t)))<0 
\] for almost every $t\in[a,b]$. Fix $q_0\in\R^{2n+1}$. The future timelike reachable set from $q_0$, denoted by $I^+(q_0)$, is the set of all points that can be reached from $q_0$ by a timelike future directed curve. Similarly, the future nonspacelike reachable set from $q_0$, denoted by $J^+(q_0)$, is the union of $\{q_0\}$ and the set of all points that can be reached from $q_0$ by a nonspacelike future directed curve. The sub-Lorentzian length of a nonspacelike future directed curve $\gamma\colon[a,b]\longrightarrow\mathbb{R}^{2n+1}$ is defined by 
\[ 
  L(\gamma) = \int_a^b \sqrt{-g(\dot{\gamma}(t),\dot{\gamma}(t))}\,dt\text{.} 
\] 
Such a curve is called a length maximizer, or a maximizing geodesic, if it has maximal length among all nonspacelike future directed curves joining $\gamma(a)$ to $\gamma(b)$. We will also use the notion of a Hamiltonian geodesic. A horizontal curve $\gamma\colon I\to\R^{2n+1}$ is called a Hamiltonian geodesic if it is the projection onto $\R^{2n+1}$ of an integral curve of the Hamiltonian vector field associated with the geodesic Hamiltonian of the Heisenberg sub-Lorentzian structure. The geodesic Hamiltonian is defined precisely in Section \ref{SecPrelimSL}. The exponential map with pole at a point $q_0$, $\exp_{q_0}\colon\mathcal{D}\longrightarrow\R^{2n+1}$, is defined by 
\[ 
 \exp_{q_0}(\lambda) = \pi\circ\Phi^1(\lambda), 
\] 
where \(\Phi^t\) is the (local) flow of the Hamiltonian vector field associated with the geodesic Hamiltonian, $\pi\colon T^*\R^{2n+1}\longrightarrow\R^{2n+1}$ is the canonical projection, and $\mathcal{D}$ is the set of all such covectors $\lambda\in T_{q_0}^*\R^{2n+1}$ that $t\longrightarrow\Phi^t(\lambda)$ is defined on $[0,1]$. Note that, as explained below, other conventions for defining $\exp_{q_0}$ are also possible. Finally, by an isometry of $(\R^{2n+1},H,g)$, we mean a diffeomorphism of $\R^{2n+1}$ onto itself that preserves both $H$ and $g$, while an infinitesimal isometry is a vector field whose (local) flow is composed of (local) isometries.

\smallskip

The Heisenberg structure on $\R^3$ was studied in \cite{Gro2006,SachkSachk}. To the best of the author’s knowledge, the case $n\geq 2$ has so far been considered only in \cite{Huang}, where the authors compute the timelike Hamiltonian geodesics, estimate the set $I^+(0)$ without determining its boundary $\partial I^+(0)$, and find the first timelike conjugate locus of the origin. The aim of the present paper is to provide a more comprehensive treatment of the case $n\geq 2$. In particular, we give a complete characterization of the class of geometrically optimal curves.

\smallskip

We begin with the computation of the isometry group.
\begin{theorem}\label{Th1}
 The group $G$ of isometries of $(\R^{2n+1},H,g)$ is 
 \begin{align}\label{IsGr}
  G = \R^{2n+1} \rtimes \big(\mathrm{O}(1,2n-1)\cap\mathrm{Sp}(2n)\big)
 \end{align}
and has dimension $n^2+3$.
\end{theorem}
After determining the dimension of $G$, we exhibit a basis for the Lie algebra of infinitesimal isometries -- see Section \ref{SecIsom}.

We then prove the following result concerning the existence of maximizing geodesics for the Heisenberg sub-Lorentzian structure.
\begin{theorem}\label{Th2}
 For every $q_0\in\R^{2n+1}$ and every $q\in J^+(q_0)\backslash\{q_0\}$, there exists a maximizing geodesic joining $q_0$ to $q$. If $q\in I^+(q_0)$, then every length maximizer joining $q_0$ to $q$ is a timelike future directed Hamiltonian geodesic and hence is a smooth curve. If $q\in\partial J^+(q_0)$, then every length maximizer joining $q_0$ to $q$ is a null future directed Hamiltonian geodesic and consequently is also smooth.
\end{theorem}
By a smooth curve in the above theorem we mean a $1$-dimensional embedded submanifold. The assertion of Theorem \ref{Th2} concerning timelike maximizers admits the following generalization.
\begin{theorem}\label{Th3}
 Let $(M,H,g)$ be a sub-Lorentzian manifold, where $H$ is a generic distribution of rank $k\geq 3$ in the sense of \cite{Trelat}. Suppose that $q_0,q\in M$, $q\in I^+(q_0)$, and that there exists a length maximizer joining $q_0$ to $q$. Then every such maximizer is a timelike future directed Hamiltonian geodesic.
\end{theorem}

Furthermore, we determine the future null and timelike conjugate and cut loci of the origin; their definitions are analogous to those in the classical setting. For every positive integer $k$, there exist numbers $\alpha_k$ and $\beta_k$ satisfying
\[
 2k\pi<\beta_k<\alpha_k<2k\pi+\tfrac{\pi}{2}
\]
such that upon defining the sets $\mathcal{A}_k$, $\mathcal{B}_k$, $\tilde{\mathcal{A}}_k$, $\tilde{\mathcal{B}}_k$ by the formulas:
\begin{align*}
 \begin{aligned} 
  &\mathcal{A}_k = \Bigl\{\Big(\tfrac{(\mathrm{sgn\,}r)\sinh 2k\pi-B_1(\cosh 2k\pi -1)}{r},\tfrac{(\mathrm{sgn\,}r)B_1\sinh 2k\pi-(\cosh 2k\pi-1)}{r},0,\dots,0,\\
  &\tfrac{1-B_1^2}{2r^2}(\mathrm{sgn\,}r)\sinh 2k\pi \Big)\in\R^{2n+1}: -1<B_1<1,\; r\neq 0 \Bigr\} \text{,}
 \end{aligned}
\end{align*}
\begin{align*}
 \begin{aligned} 
  &\mathcal{B}_k = \Bigl\{\Big(\tfrac{(\mathrm{sgn\,}r)\sinh\alpha_k - B_1(\cosh\alpha_k - 1)}{r},\tfrac{(\mathrm{sgn\,}r)B_1\sinh\alpha_k - (\cosh\alpha_k-1)}{r},0,\dots,0,\\
  &(\mathrm{sgn\,}r)\tfrac{1-B_1^2}{2r^2}\sinh\alpha_k - (\mathrm{sgn\,}r)\tfrac{\sum\nolimits_j(A_j^2+B_j^2)}{2r^2}\sin\alpha_k\Big)\in\R^{2n+1} : \\
  & -1<B_1<1,\; -1+B_1^2 + \textstyle\sum\nolimits_j(A_j^2+B_j^2)=0, \; r\neq 0 \Bigr\} \text{,}
 \end{aligned}
\end{align*}
\begin{align*}
 \begin{aligned} 
  &\tilde{\mathcal{A}}_k = \Big\{\Big(\tfrac{(\mathrm{sgn\,}r)\sinh 2k\pi-B_1(\cosh 2k\pi-1)}{r},\tfrac{(\mathrm{sgn\,}r)B_1\sinh 2k\pi - (\cosh 2k\pi-1)}{r},0,\dots,0,\\
  &(\mathrm{sgn\,}r)\tfrac{1-B_1^2}{2r^2}\sinh 2k\pi + (\mathrm{sgn\,}r)\tfrac{-1+B_1^2+\sum_j(A_j^2+B_j^2)}{r^2} k\pi\Big) 
  \in\R^{2n+1}:\\
  &B_1^2 + {\textstyle\sum\nolimits_j}(A_j^2+B_j^2)<1,\; {\textstyle\sum\nolimits_j}(A_j^2+B_j^2)\geq 0,\; r\neq 0 \Big\} \text{,}
 \end{aligned}
\end{align*}
and
\begin{align*}
 \begin{aligned} 
  &\tilde{\mathcal{B}}_k = \Big\{\Big(\tfrac{(\mathrm{sgn\,}r)\sinh\beta_k-B_1(\cosh\beta_k-1)}{r},\tfrac{(\mathrm{sgn\,}r)B_1\sinh\beta_k - (\cosh\beta_k-1)}{r},0,\dots,0,\\
  &(\mathrm{sgn\,}r)\tfrac{1-B_1^2}{2r^2}\sinh\beta_k -(\mathrm{sgn\,}r)\textstyle\sum_{j=2}^n \frac{A_j^2+B_j^2}{2r^2}\sin\beta_k 
  + (\mathrm{sgn\,}r)\tfrac{-1+B_1^2+\sum_{j=2}^n (A_j^2+B_j^2)}{2r^2}\beta_k\Big) 
  \in\R^{2n+1}: \\
  &B_1^2 + {\textstyle\sum\nolimits_j}(A_j^2+B_j^2)<1,\; {\textstyle\sum\nolimits_j}(A_j^2+B_j^2)>0,\; r\neq 0 \Big\} \text{,}
 \end{aligned}
\end{align*}
the following two theorems hold.

\begin{theorem}\label{Th4}
The future null cut and conjugate loci of the origin are, respectively,
\begin{gather*}
 \mathrm{Cut}^{null}(0)= \mathcal{A}_1\text{,}
\end{gather*}
\begin{gather*}
 \mathrm{Conj}^{null}(0)= \bigcup_{k=1}^{\infty} (\mathcal{A}_k\cup\mathcal{B}_k)\cup \{(t,\pm t,0,\dots,0):\; t>0\} \text{.}
\end{gather*}
\end{theorem}

\begin{theorem}\label{Th5}
The future timelike cut and conjugate loci of the origin are, respectively,
\begin{gather*}
 \mathrm{Cut}^{t}(0)= \tilde{\mathcal{A}}_1\text{,}
\end{gather*}
and
\begin{gather*}
 \mathrm{Conj}^{t}(0)= \bigcup_{k=1}^{\infty} (\tilde{\mathcal{A}}_k\cup \tilde{\mathcal{B}}_k)\text{.}
\end{gather*}
\end{theorem}

Finally, the preceding results allow us to give a complete description of geometrically optimal curves for the Heisenberg sub-Lorentzian metric on $\R^{2n+1}$. Recall that a nonspacelike future directed curve $\gamma:I\longrightarrow\R^{2n+1}$, where $I=[a,b]$ or $I=[a,+\infty)$, is said to be geometrically optimal if $\gamma(t)\in\partial J^+(\gamma(a))$ for every $t\in I$.

\begin{theorem}
Suppose that $\gamma$ is a geometrically optimal curve starting from the origin. Then $\gamma$ is a null future directed Hamiltonian geodesic. More precisely, up to reparametrization, one of the following alternatives holds:
\begin{enumerate}[label=\textnormal{(\arabic*)}]
    \item\label{It1} The curve \(\gamma\) is defined on \([0,+\infty)\) and is generated by an initial covector of one of the
    following forms:
    \[
    (A_1,B_1,\dots,A_n,B_n,0),
    \]
    where
    \[
    -A_1^2+B_1^2+\sum_{j=2}^n(A_j^2+B_j^2)=0, \quad A_1<0 \text{,}
    \]
    or
    \[
    (-1,\pm1,0,\dots,0,r) \text{;}
    \]

    \item The curve \(\gamma\) is defined on $[0,\tfrac{2\pi}{|r|}]$ and is generated by an initial covector
    \[
    (A_1,B_1,\dots,A_n,B_n,r)\text{,}
    \]
    where
    \[
    -A_1^2+B_1^2+\sum_{j=2}^n(A_j^2+B_j^2)=0, \quad A_1<0, \qquad |B_1|<1, \qquad r\ne0 \text{.}
    \]
\end{enumerate}
\end{theorem}
Note that in the case specified in \ref{It1}, $\gamma$ has the form given in \eqref{GeodEqsrzero} below. We also emphasize that geometrically optimal curves for the Heisenberg sub-Lorentzian metric on $\R^3$ are generally not smooth.

\section{Sub-Lorentzian structures} \label{SecPrelimSL} 
\subsection{The general case}\label{subsgen}
In this subsection we recall some basic notions of sub-Lorentzian geometry -- cf. \cite{Gro2002,Gro2009}.

Let $M$ be a smooth (i.e., of class $C^{\infty}$) $n$-dimensional manifold, and let $H$ be a smooth bracket-generating distribution of constant rank equal to $k$ on $M$. Suppose that $g$ is a Lorentzian metric on $H$. The pair $(H,g)$ is called a \textit{sub-Lorentzian metric} on $M$, and the triple $(M,H,g)$ is a \textit{sub-Lorentzian manifold}.

Fix a sub-Lorentzian manifold \((M,H,g)\). We use the standard convention from classical Lorentzian geometry. Given $q\in M$, a vector \(v\in H_q\) is called 
\begin{itemize}[label= \tiny$\bullet$]
\item \textit{timelike} (\textit{t.} for short) if \(g(v,v)<0\); 
\item \textit{spacelike} if \(g(v,v)>0\) or \(v=0\); 
\item \textit{null} if \(g(v,v)=0\) and \(v\ne0\); 
\item \textit{nonspacelike} (\textit{nspc.} for short) if \(g(v,v)\leq0\) and \(v\ne0\). 
\end{itemize}

A vector field on $M$ taking values in $H$ is called a \textit{horizontal vector field}. A continuous horizontal vector field on $M$ that is timelike at every point is called a \textit{time orientation} of $(M,H,g)$. If such a vector field, say $X$, is given, we say that $(M,H,g)$ is time-oriented by $X$. Suppose that $(M,H,g)$ is time-oriented by $X$. Then a nspc. $v\in H_q$ is called
\begin{itemize}[label= \tiny$\bullet$]
 \item \textit{future directed} (\textit{f.d.} for short) if $g(v,X(q))<0$;
 \item \textit{past directed} if $g(v,X(q))>0$.
\end{itemize}
Many features of the geometry of $(M,H,g)$ are reflected in the properties of horizontal curves. A curve $\gamma\colon I\longrightarrow M$ is called \textit{horizontal} if it is absolutely continuous and $\dot{\gamma}(t)\in H_{\gamma(t)}$ for almost every $t\in I$; here $I$ may be an arbitrary interval. A horizontal curve $\gamma\colon I\longrightarrow M$ is called \textit{timelike} (respectively, \textit{timelike future directed}, \textit{nonspacelike}, etc.) if $\dot{\gamma}(t)$ is so for almost every $t$ in $I$.\smallskip

From now on, unless otherwise stated, all curves are assumed to be horizontal. Moreover, in accordance with the abbreviations introduced above, we write \textit{t.f.d.}, \textit{nspc.f.d.}, and \textit{null f.d.} for \textit{timelike future directed}, \textit{nonspacelike future directed}, and \textit{null future directed}, respectively. Thus, for instance, by a t.f.d. curve we mean a horizontal curve which is timelike future directed. We also adopt the following convention: by a smooth curve we mean a horizontal curve $\gamma\colon I\longrightarrow M$ of class $C^{\infty}$, where $I$ is an open (resp. closed) interval and $\gamma(I)$ is an embedded submanifold (resp. an embedded submanifold with boundary).\smallskip

We now define a family of sets describing the causal structure of $(M,H,g)$. Fix a point $q_0\in M$. We denote by $I^+(q_0)$ the set of all points in $M$ that can be reached from $q_0$  by a t.f.d. curve. $I^+(q_0)$ is called the \textit{future timelike reachable set from $q_0$}. Similarly, $J^+(q_0)$ (resp. $N^+(q_0)$) is defined as the union of $\{q_0\}$ and the set of all points in $M$ that can be reached from $q_0$ by an nspc.f.d. (resp. by a null f.d.) curve. $J^+(q_0)$ (resp. $N^+(q_0)$) is called the \textit{future nonspacelike reachable set from $q_0$} (resp. \textit{future null reachable set from $q_0$}). There are also local variants of the definitions of the reachable sets, but we will not need them in this paper. For general properties of reachable sets the reader is referred to \cite{Gro2009}.

Let $\gamma\colon [a,b]\longrightarrow M$ be an nspc.f.d. curve. We define its \textit{sub-Lorentzian length}, $L(\gamma)$, by the formula
\[
 L(\gamma) = \int_a^b|g(\dot{\gamma}(t),\dot{\gamma}(t))|^{\frac12}dt\text{.}
\]
An nspc.f.d. curve $\gamma\colon [a,b]\longrightarrow M$ is called a \textit{maximizing geodesic}, a \textit{length maximizer}, or simply a \textit{maximizer} if it is longest among all nspc.f.d. curves joining $\gamma(a)$ and $\gamma(b)$. By a \textit{geodesic}, we mean an nspc.f.d. curve $\gamma\colon I\longrightarrow M$ such that, for every sufficiently short subinterval $I_1\subset I$, the restriction $\gamma_{|I_1}$ is a maximizer.

Let $X_1,\dots,X_k$ be an orthonormal frame for $H$, with $X_1$ timelike, defined in some open set $U$ in $M$. As is readily seen, the function $\mathcal{H}:T^*U\longrightarrow\R$ defined by
\[
 \mathcal{H}(q,\lambda)= -\frac12\langle\lambda,X_1(q) \rangle^2 + \frac12\sum_{i=2}^k\langle\lambda,X_i(q) \rangle^2
\]
is independent of the choice of local orthonormal frame. Therefore, $\mathcal{H}$ is well defined on the entire cotangent bundle $T^*M$, and it is called the \textit{geodesic Hamiltonian} induced by the structure $(H,g)$. Let $\overrightarrow{\mathcal{H}}$ denote the Hamiltonian vector field associated with $\mathcal{H}$ and let $\Phi^t$ stand for its (local) flow. Every curve of the form $\gamma(t) = \pi \circ \Phi^t(q_0,\lambda_0)$, where $\pi\colon T^*M\longrightarrow M$ is the canonical projection, is called a \textit{Hamiltonian geodesic}. Every Hamiltonian geodesic is a smooth curve and has constant causal character. The latter property means that if $\gamma\colon I\longrightarrow M$ is such a curve and the vector $\dot{\gamma}(t_0)$ is, for instance, timelike for some $t_0\in I$, then $\dot{\gamma}(t)$ is timelike for every $t\in I$. It can be also proved that every nspc.f.d. Hamiltonian geodesic is a geodesic in the above sense, i.e., it is locally length maximizing; see \cite{Gro2002,Gro2009}. Fix a point $q_0\in M$ and consider the set 
\[
 \mathcal{D}_{q_0} = \{\lambda\in T_{q_0}^*M:\; \text{the curve }t\longrightarrow\Phi^t(q_0,\lambda) \text{ is defined on }[0,1]\} \text{.}
\]
The mapping
\begin{align*}
 \exp_{q_0}\colon\mathcal{D}_{q_0}\longrightarrow M, \; \exp_{q_0}(\lambda) = \pi\circ\Phi^1(q_0,\lambda)
\end{align*}
is called the \textit{exponential map with pole at $q_0$}. This is the standard definition, although other equivalent approaches are possible. In this definition, the initial covectors are not normalized, whereas the final time is fixed. Below, we adopt a different convention: the initial covectors are normalized, while the final time is allowed to vary. By the \textit{future timelike conjugate locus} of a point $q_0$ we mean the set $\mathrm{Conj}^t(q_0)$ defined by
\[
 \{\lambda\in\mathcal{D}_{q_0}:\; \mathcal{H}(\lambda)<0,\,\langle \lambda,X\rangle<0, \, \text{and } d_{\lambda}\exp_{q_0} \text{ is singular} \} \text{.}
\]
Similarly, the \textit{future null conjugate locus} of $q_0$, $\mathrm{Conj}^{null}(q_0)$, is defined as
\[
 \{\lambda\in\mathcal{D}_{q_0}:\; \mathcal{H}(\lambda)=0,\,\langle \lambda,X\rangle<0, \, \text{and } d_{\lambda}\exp_{q_0} \text{ is singular} \} \text{.}
\]
We next define the cut locus of $q_0$. Consider a t.f.d. geodesic $\gamma\colon[0,T_{\max})\longrightarrow M$ starting from $q_0$ and defined on its maximal interval of existence. Define its \textit{cut time}  by
\[
 t_{\gamma} = \sup\{t\in(0,T_{\max}):\; \gamma_{|[0,t]}\text{ is length maximizing}\}\text{.}
\]
The \textit{future timelike cut locus} of $q_0$ is the set of all points of the form $\gamma(t_{\gamma})$, where $\gamma$ ranges over all t.f.d. geodesics as above for which $t_\gamma<T_{\max}$. 
The future null cut locus of $q_0$ is defined analogously, with null f.d. geodesics in place of t.f.d. ones.

Suppose that $\varphi\colon U\longrightarrow \R$ is a smooth function defined on an open set $U\subset M$. By the \textit{horizontal gradient} of $\varphi$ we mean the horizontal vector field, denoted by $\nabla_H\varphi$, such that
\[
 \partial_v\varphi = g(\nabla_H\varphi(q),v)
\]
for every $q\in U$ and $v\in H_q$.
It can be shown that if $\nabla_H\varphi$ is nspc.f.d. on $U$ and $g(\nabla_H\varphi,\nabla_H\varphi)$ is constant on $U$, then the trajectories of $\nabla_H\varphi$ are unique maximizing geodesics.

A diffeomorphism $f\colon M\longrightarrow  M$ is called an \textit{isometry} of $(M,H,g)$ if $df$ preserves $H$ and the restriction $d_qf_{|H_q}\colon H_q\longrightarrow H_{f(q)}$ is a linear isometry for every $q\in M$. By an \textit{infinitesimal isometry} of $M$ we mean a vector field on $M$ whose (local) flow consists of (local) isometries of $M$. It is easy to see that time-orientation-preserving isometries map geodesics, Hamiltonian geodesics, maximizers to geodesics, Hamiltonian geodesics, maximizers, respectively. They also preserve reachable sets. For instance, if $f$ is such an isometry and $q_0\in M$, then $f(J^+(q_0)) = J^+(f(q_0))$, $f(I^+(q_0)) = I^+(f(q_0))$. 

Finally, we would like to draw the reader's attention to the notions of abnormal curves and Goh curves; see \cite{MongAb} and \cite{AgrSach,Trelat}, respectively. We will not state their definitions explicitly here; we merely note that nonconstant curves of either type do not exist in the contact case.
\smallskip

\subsection{The contact case}\label{subscont}
We now define the class of contact sub-Lorentzian manifolds. Let $M$ be a smooth $(2n+1)$-dimensional manifold, let $H$ be a smooth contact distribution on $M$, and let $g$ be a smooth Lorentzian metric on $H$. Then the pair $(H,g)$ is called a \textit{contact sub-Lorentzian structure} on $M$, whereas $(M,H,g)$ is referred to as a \textit{contact sub-Lorentzian manifold}. Contact sub-Lorentzian manifolds are, in many respects, the simplest sub-Lorentzian manifolds.

Fix a time-oriented contact sub-Lorentzian manifold $(M,H,g)$. We will describe the maximizing geodesics. The standard approach is to apply Pontryagin's maximum principle, PMP for short; see, for instance, \cite{AgrSach, Bress, Pon}. Suppose that $q,q_0\in M$ are such that $q\in I^+(q_0)$, and that there exists a maximizing geodesic $\gamma$ joining $q_0$ to $q$. It is clear that $\gamma$ cannot be a null curve. We know that all maximizers are PMP extremals and that such extremals can be either normal or abnormal. Since timelike abnormal extremals are abnormal curves \cite{Gro2009}, and nonconstant abnormal curves do not exist in the contact case, it follows that $\gamma$ must be a normal extremal.

We thus investigate normal extremals. To this end, take an orthonormal frame $X_1,\dots,X_{2n}$ defined on an open set $U\subset M$, where $X_1$ is the restriction to $U$ of a time orientation of $(M,H,g)$. Every nspc.f.d. curve in $U$ is a solution of the control system
\[
 \dot{q} = \sum_{i=1}^{2n}u_i(t)X_i\text{,} 
\]
where the set of control parameters is equal to
\[
 C = \Big\{(u_1,\dots,u_{2n}):\; -u_1^2 + \sum_{j=2}^{2n}u_j^2\leq 0,\,u_1>0\Big\} \text{,} 
\]
and the admissible controls are of class $L^{\infty}$.
The cost functional to be maximized is
\[
 \int_0^T\sqrt{u_1^2(t) - \sum_{j=2}^{2n}u_j^2(t)}dt\text{,}
\]
with final time $T$ free.
The PMP Hamiltonian $h_u\colon T^*U\longrightarrow\R$ in the normal case is given by 
\begin{align*}
 h_u(q,\lambda) = \sum_{i=1}^{2n}u_i\langle \lambda,X_i(q)\rangle + \sqrt{u_1^2 - \sum_{j=2}^{2n}u_j^2} \text{.}
\end{align*}
Since length does not depend on the choice of parameterization, we may assume that $u_1=1$. Therefore, the PMP Hamiltonian takes the form
\begin{align}\label{hamPmP}
 h_u(q,\lambda) = \langle \lambda,X_1(q)\rangle + \sum_{j=2}^{2n}u_j\langle \lambda,X_j(q)\rangle + \sqrt{1 - \sum_{j=2}^{2n}u_j^2} \text{,}
\end{align}
where the set of control parameters is now
\begin{align*}
 C_1 = \Big\{(u_2,\dots,u_{2n}):\; \sum_{j=2}^{2n}u_j^2\leq 1\Big\} \text{.} 
\end{align*}
It is not difficult to verify that the Hamiltonian \eqref{hamPmP} does not attain its maximal value on $\partial C_1$. 
Indeed, let $w=(u_2,\dots,u_{2n})\in\partial C_1$ and let $w_{\varepsilon}=\sqrt{1-\varepsilon}w$, $\varepsilon\in(0,1)$. Then $w_{\varepsilon}\in\mathrm{int}\,C_1$ and 
\[
 h_{w_{\varepsilon}}(q,\lambda) - h_w(q,\lambda) = -\frac12\varepsilon\sum_{j=2}^{2n}w_j\langle \lambda,X_j(q)\rangle + \sqrt{\varepsilon}+o(\varepsilon)>0
\]
for sufficiently small $\varepsilon>0$. Consequently, the normal extremals are timelike.
We are now in a position to state the following proposition.
\begin{proposition}\label{NormareT}
 On every sub-Lorentzian manifold, all normal extremals are t.f.d. Hamiltonian geodesics. Accordingly, in the contact case, for any two points $q_0$ and $q$ such that $q\in I^+(q_0)$, every maximizer joining $q_0$ to $q$ is a t.f.d. Hamiltonian geodesic.
\end{proposition}
\begin{proof}
 The proof that normal extremals are t.f.d. Hamiltonian geodesics closely follows that of the corresponding result in sub-Riemannian geometry. We therefore omit the proof.
\end{proof}

Let us turn to null maximizers. For a given curve $\gamma$, saying that $\gamma$ is a null maximizer is equivalent to saying that it is geometrically optimal. Let us clarify the meaning of the latter notion. An nscp.f.d. curve $\gamma\colon [a,b]\longrightarrow M$ is said to be \textit{geometrically optimal} if $\gamma([a,b])\subset\partial J^+(\gamma(a))$. Note that it is clear from control theory that if $\gamma(b)\in\partial J^+(\gamma(a))$, then necessarily $\gamma$ is geometrically optimal. We have the following theorem
\begin{theorem}[\cite{Gro2009}]\label{ThNUll}
 Every sufficiently short subarc of a null f.d. Hamiltonian geodesic is geometrically optimal. Moreover, in the contact case, every geometrically optimal curve is a concatenation of finitely many null f.d. Hamiltonian geodesics.
\end{theorem}
\begin{corollary}
 In the contact case, the sets $I^+(q)$, $q\in M$, are open.
\end{corollary}

At the end of this section let us note that the claim of Proposition \ref{NormareT} can be generalized to sub-Lorentzian manifolds $(M,H,g)$, where $H$ is a rank $k\geq 3$ generic distribution in the sense of \cite{Trelat}. Indeed, using the similar argument as in \cite[Lemma 6.3]{Gro2009}, one can prove that every t. maximizer which is an abnormal curve is a Goh curve. However, by \cite{Trelat} we know that Goh curves do not exist for this class of distribution. This proves Theorem \ref{Th3}.

\section{Geodesics of the Heisenberg sub-Lorentzian structure on $\R^{2n+1}$}\label{secHeisIsom}
As above we consider the space $\R^{2n+1}$ with coordinates $(x_1,y_1,\dots,x_n,y_n,z)$. We define the vector fields $X_i$,$Y_i$, the contact distribution $H$, and the Lorentzian metric $g$ on $H$, respectively, by
\begin{align}
 X_i = \frac{\partial}{\partial x_i} + \frac12 y_i\frac{\partial}{\partial z}, \; Y_i = \frac{\partial}{\partial y_i} - \frac12 x_i\frac{\partial}{\partial z}\text{,}
\end{align}
$i=1,\dots,n$,
\[
 H = \Sp\{X_1,Y_1,\dots,X_n,Y_n\}\text{,}
\]
and
\[
 g(X_1,X_1)=-1,\; g(Y_1,Y_1)=1,\; g(X_j,X_j)=g(Y_j,Y_j)=1,\; j=2,\dots,n\text{,}
\]
\[
 g(X_i,X_j)=g(Y_i,Y_j)=0,\; 1\leq i<j\leq n,\quad g(X_i,Y_j)=0,\; i,j=1,\dots,n \text{.}
\]
We take $X_1$ to be a time orientation. Recall that the pair $(H,g)$ is called the Heisenberg sub-Lorentzian structure on $\R^{2n+1}$.
Let us notice that $H = \ker\omega$, where
\begin{align}\label{ContForm}
 \omega = dz - \frac12\sum_{i=1}^n(y_idx_i - x_idy_i)\text{.}
\end{align}
It follows that every horizontal curve in $\R^{2n+1}$ is uniquely determined by its projection onto the $(x_1,y_1,\dots,x_n,y_n)$-plane, identified with $\R^{2n}$, up to the choice of the initial value of its $z$-coordinate. Indeed, if $\tilde{\gamma}=(x_1,y_1,\dots,x_n,y_n)\colon[a,b]\longrightarrow\R^{2n}$ is an absolutely continuous curve, then it can be lifted to a horizontal curve $\gamma=(x_1,y_1,\dots,x_n,y_n,z)\colon[a,b]\longrightarrow\R^{2n+1}$ by imposing the \textit{horizontality condition} 
\begin{align}\label{HorCon}
 z(t) = z_0 + \frac12\int_a^t \sum_{i=1}^n(y_i\dot{x}_i - x_i\dot{y}_i)\mathrm{d}\tau\text{.}
\end{align}
This fact will be used repeatedly below.

We will compute nspc.f.d. geodesics on $(\R^{2n+1},H,g)$. The geodesic Hamiltonian associated to the structure $(H,g)$ has the form
\[
 \mathcal{H} = -\frac12 f_1^2 + \frac12 g_1^2 + \frac12 \sum_{j=2}^n (f_j^2 + g_j^2),
\]
where $f_i,g_i:T^*\R^{2n+1}\longrightarrow\R$ are given by
\begin{gather*}
 f_i = \langle\lambda,X_i\rangle = p_i + \frac12y_ir\text{,} \\
 g_i = \langle\lambda,Y_i\rangle = q_i - \frac12x_ir \text{,}
\end{gather*}
$i=1,\dots,n$. Here $\lambda = (p_1,q_1,\dots,p_n,q_n,r)$ and $(x_1,y_1,\dots,x_n,y_n,z,p_1,q_1,\dots,p_n,q_n,r)$ are Darboux coordinates on $T^*\R^{2n+1}$. Observing that
\begin{align}
 \frac{\partial f_i}{\partial x_j} = \frac{\partial g_i}{\partial y_j} = 0,\; \frac{\partial f_i}{\partial y_j} = \frac12 r\delta^i_j, \; \frac{\partial g_i}{\partial x_j} = -\frac12 r\delta^i_j\text{,}
\end{align}
$i,j=1\dots,n$, we can explicitly write out Hamiltonian equations and obtain

\begin{align}\label{geq1}
 \left\{\begin{array}{l}
         \dot{f}_1 = rg_1 \\
         \dot{g}_1 = rf_1
        \end{array} \right. \text{,}\;\;
        \left\{\begin{array}{l}
         \dot{f}_j = rg_j \\
         \dot{g}_j = -rf_j
        \end{array} \right.\text{,}
\end{align}
and
\begin{align}\label{geq2}
 \left\{\begin{array}{l}
         \dot{x}_1 = -f_1 \\
         \dot{y}_1 = g_1
        \end{array} \right. \text{,}\;\;
         \left\{\begin{array}{l}
         \dot{x}_j = f_j \\
         \dot{y}_j = g_j
        \end{array} \right.\text{,}
\end{align}
$j=2,\dots,n$. As will soon become clear, in order to investigate the behavior of geodesics, it suffices to compute those starting from the origin. Thus, we consider geodesics with initial conditions $x_i(0)=0$, $y_i(0)=0$, $p_i(0)=A_i$, $q_i(0)=B_i$, $i=1,\dots,n$, and $z(0)=0$, $r(0)=r$. After calculations we find
\begin{align}\label{geq3}
 \left\{\begin{array}{l}
         f_1(t) = A_1\cosh(rt) + B_1\sinh(rt) \\
         g_1(t) = B_1\cosh(rt) + A_1\sinh(rt) 
        \end{array} \right. \text{,}\;\;
         \left\{\begin{array}{l}
         f_j(t) = A_j\cos(rt) + B_j\sin(rt) \\
         g_j(t) = B_j\cos(rt) - A_j\sin(rt)
        \end{array} \right.\text{,}
\end{align}
$j=2,\dots,n$. Note that since $\dot{r} = \frac{\partial\mathcal{H}}{\partial z}=0$, $r$ is always constant. Applying the horizontality condition \eqref{HorCon} we obtain the equation for the $z$-coordinate of geodesics:
\begin{align}\label{geq4}
 \dot{z} = \tfrac{-A_1^2+B_1^2+\sum_{j=2}^n(A_j^2+B_j^2)}{2r} + \tfrac{A_1^2-B_1^2}{2r}\cosh(rt) - \sum_{j=2}^n\tfrac{A_j^2+B_j^2}{2r}\cos(rt)\text{.}
\end{align}
Using (\ref{geq1}), (\ref{geq2}), (\ref{geq3}), (\ref{geq4}) we finally obtain
\begin{align}\label{GeodEqsrzero}
\left\{\begin{array}{l}
 x_1(t) = - A_1 t \smallskip \\
 y_1(t) = B_1 t \smallskip \\  
 x_j(t) = A_j t \smallskip \\
 y_j(t) = B_j t \smallskip \\
 z(t) = 0
\end{array}\right.
\end{align}
for $r=0$, and 
\begin{align}\label{GeodEqsrnotzero}
\left\{\begin{array}{l}
 x_1(t) = -\frac{A_1}{r}\sinh(rt) - \frac{B_1}{r}(\cosh(rt)- 1) \bigskip \\
 y_1(t) = \frac{B_1}{r}\sinh(rt) + \frac{A_1}{r}(\cosh(rt) -1)  \bigskip \\  
 x_j(t) = \frac{A_j}{r}\sin(rt) - \frac{B_j}{r}(\cos(rt) -1) \bigskip \\
 y_j(t) = \frac{B_j}{r}\sin(rt) + \frac{A_j}{r}(\cos(rt) -1) \bigskip \\
 z(t) = \frac{A_1^2-B_1^2}{2r^2}\sinh(rt) - \sum_{j=2}^n \frac{A_j^2+B_j^2}{2r^2}\sin(rt) + \frac{-A_1^2+B_1^2+\sum_{j=2}^n (A_j^2+B_j^2)}{2r}t
\end{array}\right.
\end{align}
for $r\ne 0$, where again $j=2,\dots,n$. In the equations above we abbreviate
\[
 x_1(t;A_1,B_1,\dots,A_n,B_n,r) \text{,} 
\]
where the initial covector is explicitly indicated, as $x_1(t)$. We use the same convention for the remaining coordinates. Let us remark that the equations \eqref{GeodEqsrzero} can be obtained from those in \eqref{GeodEqsrnotzero} by letting $r$ tend to zero. Moreover, the formulas in \eqref{GeodEqsrzero} and \eqref{GeodEqsrnotzero} together define functions that depend smoothly on all their arguments, including $r=0$. 

Below we gather some useful facts about geodesics (starting from the origin) which easily follow from the Hamiltonian equations.
\begin{remark}\label{Rem1}
\noindent
\begin{enumerate}
 \item \label{a1} A nspc.f.d. geodesic $\gamma$, with $\gamma(0)=0$, corresponds to an initial covector $\lambda(0) =(A_1,B_1,\dots,A_n,B_n,r)\in T_0^*\R^{2n+1}$ satisfying $-A_1^2 + B_1^2 + \sum_{j=2}^n (A_j^2+B_j^2)\leq 0$, $-A_1>0$. Its initial velocity is $\dot{\gamma}(0)=-A_1X_1(0) + B_1Y_1(0) + \sum_{j=2}^n(A_jX_j(0)+B_jY_j(0))$. 
  
 \item  Suppose that $(\gamma(t),\lambda(t))$ and $(\tilde{\gamma}(t),\tilde{\lambda}(t))$ are solutions of the equation 
 \[
 (\gamma(t),\lambda(t))^{\boldsymbol{\cdot}} = \overrightarrow{\mathcal{H}}(\gamma(t),\lambda(t)) 
 \]
with initial conditions $(\gamma(0),\lambda(0))$ and $(\gamma(0),s\lambda(0))$, respectively, $s\in\R$. Then $\tilde{\gamma}(t) = \gamma(st)$ and $\tilde{\lambda}(t) = s\lambda(st)$. Thus replacing $\lambda(0)$ with $s\lambda(0)$ results in the same unparameterized geodesic, but with a different parameterization.
 
 \item The family of initial covectors $(A_1,B_1,\dots,A_n,B_n,sr)\in T_0^*\R^{2n+1}$, $s\in\R$, produces different geodesics with the same initial velocity at zero as in \eqref{a1}.
\end{enumerate}

\end{remark}

A closer look at the above formulas leads to the following observations.
\begin{proposition}
 If $j\in\{2,\dots,n\}$, then the projections of the geodesics on the $(x_j,y_j)$-plane are either straight lines when $r=0$, or circles (possibly degenerate to a point) given by the equation
 \begin{align*}
  \Big(x_j-\frac{B_j}{r}\Big)^2 + \Big(y_j+\frac{A_j}{r}\Big)^2 = \frac{A_j^2+B_j^2}{r^2}
 \end{align*}
in case $r\ne 0$.
\end{proposition}

\begin{proposition}
 The projections of the nonspacelike geodesics on the $(x_1,y_1)$-plane are either straight lines when $r=0$, or branches of hyperbolas (possibly degenerate to a pair of perpendicular straight lines) given by the equation
 \begin{align*}
  \Big(x_1 - \frac{B_1}{r}\Big)^2 - \Big(y_1 + \frac{A_1}{r}\Big)^2 = \frac{-A_1^2 + B_1^2}{r^2}
 \end{align*}
in case $r\ne 0$.
\end{proposition}
In the above propositions the degenerate situation arises when $A_1=\pm B_1$. Then, necessarily, $A_j=B_j=0$, $j=2,\dots,n$.

\section{Isometries of $(\R^{2n+1},H,g)$}\label{SecIsom}
In this section we prove Theorem \ref{Th1} and show how to obtain explicit formulas for isometries. Note that the formula \eqref{IsGr} in Theorem \ref{Th1} follows from \cite[Subsection 3.2]{GrKrESAIM}. We will compute the dimension of $G$. It is enough to prove the following proposition.
\begin{proposition}
 \[
  \dim \big(\mathrm{O}(1,2n-1)\cap \mathrm{Sp}(2n)\big) = n^2 - 2n + 2 \text{.}
 \]
\end{proposition}
\begin{proof}
We will work at the level of Lie algebras. Let
 \[
  J = \begin{pmatrix}
       \large 0 & -I_n \\
       I_n & \large 0
      \end{pmatrix},\;
 \]
where $I_n$ denotes the $n\times n$ identity matrix, and let $Q = \mathrm{diag}(-1,1,\dots,1)$ be the $2n\times 2n$ diagonal matrix with diagonal entries $-1,1,\dots,1$. Suppose that $M\in \mathfrak{o}(1,2n-1)\cap \mathfrak{sp}(2n)$. Then, by definition, 
\[ 
  M^TQ+QM=0 \qquad\text{and}\qquad M^TJ+JM=0 \text{.} 
\]
It follows that $M=(m_{ij})_{i,j=1,\dots,2n }$ can be written simultaneously in the following two block forms:
\begin{align}\label{MatAlg}
 M = \begin{pmatrix}
      0 & X^T \\
      X & D
     \end{pmatrix} = \;
     \begin{pmatrix}
      A & B \\
      C & -A^T
     \end{pmatrix}\text{,}
\end{align}
where $0\in\R$, $D\in\mathfrak{o}(2n-1)$, $X\in\R^{2n-1}$, $A,B,C$ are $n\times n$ matrices, and $B^T=B$, $C^T=C$. First, notice that the only entry in the first row of $M$ that may be nonzero is $m_{1,n+1}$. Indeed, the symmetry relations resulting from the two representations in \eqref{MatAlg} give
\begin{align}\label{MatSymm}
 \begin{array}{c}
  m_{1i} = m_{i1}, \quad i=1,\dots,2n  \text{,} \smallskip\\
  m_{1,j} = - m_{n+j,n+1}, \quad j=2,\dots,n \text{,}\smallskip \\
 m_{1,n+j} = m_{j,n+1}, \quad j=2,\dots,n \text{,}\smallskip \\
 m_{n+1,j} = m_{n+j,1}, \quad j=2,\dots,n \text{,}\smallskip \\
 m_{rs} = -m_{sr}, \quad r,s = 2,\dots,2n \text{.}
 \end{array}
\end{align}
Consider, for instance, $m_{1,n+j}$, where $j=2,\dots,n$. We have
\[
 m_{1,n+j} = m_{j,n+1}
\]
and 
\[
 m_{1,n+j} = m_{n+j,1} = m_{n+1,j} = -m_{j,n+1}\text{,}
\]
from which it follows that $m_{1,n+j}=0$. Thus $X$ is described by $1$ parameter.

Next, after deleting the first row and the first column of $A$, which have already been determined by \eqref{MatSymm}, we are left with an $n-1 \times n-1$ skew-symmetric matrix, which is described by $\frac{(n-1)(n-2)}{2}$ parameters. Similarly, after deleting the first row and the first column of $B$, which have also been already determined by \eqref{MatSymm}, we are left with an $n-1 \times n-1$ symmetric matrix, which is described by $\frac{(n-1)n}{2}$ parameters. Finally, using the mentioned symmetry relations, we see that $C$ is determined by $B$. Consequently, the dimension of $\mathfrak{o}(1,2n-1)\cap \mathfrak{sp}(2n)$ equals
\begin{gather*}
  1 + \frac{(n-1)n}{2} + \frac{(n-1)(n-2)}{2} = n^2 - 2n + 2\text{,}
\end{gather*}
 which completes the proof.
\end{proof}

One method to explicitly compute the group of isometries is to apply the procedure described in \cite[Subsection 3.2]{GrKrESAIM}. The reader should note, however, that the ordering of the coordinates in \cite{GrKrESAIM} differs from that used in the present paper. Alternatively, one can compute the generators of the Lie algebra of infinitesimal isometries and then, if necessary, integrate them to obtain the formulas for time-orientation-preserving isometries. We adopt the latter approach. Following e.g. \cite[Proposition 2.2]{Gro2025} (the formulation in the indefinite case is exactly the same), a vector field $Z=\sum_{i=1}^n(f_iX_i+g_iY_i)+h\frac{\partial}{\partial z}$, where the functions $f_i,g_i,h$ are smooth, is an infinitesimal isometry of $(\R^{2n+1},H,g)$ if and only if it satisfies the following system of $2n^2+3n$ PDE's:
\begin{align}\label{EqInfIsom}
\begin{array}{c}
 g([Z,X_i],X_j) + g(X_i,[Z,X_j]) = 0, \text{,} \smallskip \\
 g([Z,X_r],Y_s) + g(X_r,[Z,Y_s]) = 0 \text{,} \smallskip \\
 g([Z,Y_i],Y_j) + g(Y_i,[Z,Y_j]) = 0 \text{,} \smallskip \\
 \omega([Z,X_s]) = \omega([Z,Y_s]) = 0\text{,}
\end{array}
\end{align}
where $1\leq r,s\leq n$, $1\leq i\leq j \leq n$ and $\omega$ is defined in \eqref{ContForm}. Solving \eqref{EqInfIsom}, which is straightforward in the Heisenberg case, yields the following list of the mentioned generators:
\begin{gather*} 
 \frac{\partial}{\partial x_i} - \frac12 y_i\frac{\partial}{\partial z},\; i=1,\dots,n \text{,}\\
 \frac{\partial}{\partial y_i} + \frac12 x_i\frac{\partial}{\partial z},\; i=1,\dots,n \text{,}\\
 \frac{\partial}{\partial z} \text{,}\\
 y_1\frac{\partial}{\partial x_1} + x_1\frac{\partial}{\partial y_1} \text{,}\\
 -y_j\frac{\partial}{\partial x_j} + x_j\frac{\partial}{\partial y_j},\; j=2,\dots,n \text{,}\\
 -y_j\frac{\partial}{\partial x_k} + x_k\frac{\partial}{\partial y_j} - y_k\frac{\partial}{\partial x_j} + x_j\frac{\partial}{\partial y_k} ,\; j,k=2,\dots,n,\, k<j \text{,}\\
 -x_j\frac{\partial}{\partial x_k} + x_k\frac{\partial}{\partial x_j} - y_j\frac{\partial}{\partial y_k} + y_k\frac{\partial}{\partial y_j} ,\; j,k=2,\dots,n,\, k<j \text{.}
\end{gather*}
By integrating the first $2n+1$ fields in the above list, we see that the group of time-orientation-preserving isometries acts transitively on $\R^{2n+1}$. We will use this fact in what follows.

We conclude this section by noting that $(\mathbb{R}^{2n+1},H,g)$ admits a Lie group interpretation; see, e.g., \cite{GrKrESAIM}. Under this interpretation, the first $2n+1$ vector fields in the above list correspond to left translations.

\section{Reachable set}
\subsection{Preliminary remarks}\label{PrelRem}
Let us fix a positive integer $n\geq 2$.
We are interested in the reachable sets $I^+(q_0)$, $J^+(q_0)$ for a given $q_0\in\R^{2n+1}$. Since, as it was mentioned above, the group of time-orientation-preserving isometries acts transitively on $\R^{2n+1}$, it is enough to concentrate on the sets $I^+(0)$ and $J^+(0)$. 

Let $\pi:\R^{2n+1}\longrightarrow\R^{2n}$ be the projection 
\[
 (x_1,y_1,\dots,x_n,y_n,z)\longrightarrow(x_1,y_1,\dots,x_n,y_n,0)
\]
regarded as a mapping from $(\R^{2n+1},H,g)$ onto $(\R^{2n},g_0)$, where $g_0$ is the standard Minkowski scalar product on $\R^{2n}$. Here and below we identify $\R^{2n}$ with $\{z=0\}\subset\R^{2n+1}$.
Let us take $\frac{\partial}{\partial x_1}$ as a time orientation of $(\R^{2n},g_0)$. Since $d\pi(X_i) = \frac{\partial}{\partial x_i}$, $d\pi(Y_i) = \frac{\partial}{\partial y_i}$, the restriction
\[
 d_q\pi:H_q\longrightarrow T_{\pi(q)}\R^{2n}
\]
is a time-orientation-preserving linear isometry for every point $q\in\R^{2n+1}$. It follows that $\pi$ maps t.f.d. (resp. null f.d.) curves in $(\R^{2n+1},H,g)$ onto t.f.d. (resp. null f.d.) curves in $(\R^{2n},g_0)$. On the other hand, every horizontal curve in $\R^{2n+1}$ is determined by its projection onto $\R^{2n}$ modulo its starting point --- cf. Section \ref{secHeisIsom} and formula \eqref{HorCon}.  
This observation yields the following remark.
\begin{remark}\label{rem2}
 $\pi$ defines a bijection between the set of t.f.d. (resp., nspc.f.d) curves in $(\R^{2n+1},H,g)$ starting from the origin and the set of t.f.d. (resp., nspc.f.d) curves in $(\R^{2n},g_0)$ starting from the origin. In particular, 
 \begin{align}\label{min1}
 \pi(J^+(0)) = \{(x_1,y_1,\dots,x_n,y_n,0): -x_1^2+y_1^2+\sum_{j=2}^n(x_j^2+y_j^2)\leq 0,\, x_1\geq 0\}\text{,} 
 \end{align}
and
\begin{align}\label{min2}
 \pi(I^+(0)) = \{(x_1,y_1,\dots,x_n,y_n,0): -x_1^2+y_1^2+\sum_{j=2}^n(x_j^2+y_j^2)< 0,\, x_1> 0\}
\end{align}
\end{remark}
\begin{proof}
 To see the inclusions "$\subset$", it is enough to recall that the right-hand side of (\ref{min1}) (resp. (\ref{min2})) is the causal future (resp. chronological future) of the origin for the Minkowski metric on $\R^{2n}$. The fact that equality holds in both (\ref{min1}) and (\ref{min2}) follows from (\ref{GeodEqsrzero}) and Remark \ref{Rem1}.
\end{proof}
By Remark \ref{rem2} we also have
\[
 \partial J^+(0)\cap\{z=0\} = \{(x_1,y_1,\dots,x_n,y_n,0): -x_1^2+y_1^2+\sum_{j=2}^n(x_j^2+y_j^2) = 0,\, x_1\geq 0\} \text{.}
\]

Before proceeding further let us make the following easy but useful observations (cf. Subsection \ref{subscont}). Let $\Phi\colon\R^{2n+1}\longrightarrow\R$ be defined by $\Phi(x_1,y_1,\dots,x_n,y_n,z)=x_1$. Suppose that $\gamma:I\longrightarrow\R^{2n+1}$ is an arbitrary nspc.f.d. curve defined on an interval $I$, i.e.  
\[
 \dot{\gamma}(t) = \sum_{i=1}^n (u_i(t)X_i(\gamma(t)) + v_i(t)Y_i(\gamma(t))) \text{,}
\]
where $u_i,v_i$ are of class $L^{\infty}$ and such that 
\[
 -u_1(t)^2+v_1(t)^2 + \sum_{j=2}^n (u_j(t)^2 + v_j(t)^2)\leq 0, \; u_1(t)>0\text{,}
\]
a.e. Since $\nabla_H\Phi$ is timelike past directed, it is clear that the function $t\longrightarrow\Phi(\gamma(t))$ is increasing. It follows that every such $\gamma$ can be reparameterized so as to have $\Phi(\gamma(t)) \equiv t$ (we keep the same notation for the reparameterization), that is,  so as to have
\begin{align}\label{AffCon}
 \dot{\gamma}(t) = X_1(\gamma(t)) + \tilde{v}_1(t)Y_1(\gamma(t)) + \sum_{j=2}^n (\tilde{u}_j(t)X_j(\gamma(t)) + \tilde{v}_j(t)Y_j(\gamma(t))) \text{,}
\end{align}
where $v_1(t)^2 + \sum_{j=2}^n (u_j(t)^2 + v_j(t)^2)\leq 1$. In this way, we can reduce the study of nspc.f.d. curves to the study of trajectories of the corresponding control affine system with controls taking values in the unit ball. 

\subsection{Existence of maximizers}
Fix a point $\bar{q}=(\bar{x}_1,\bar{y}_1,\dots,\bar{x}_n,\bar{y}_n,\bar{z})\in J^+(0)$ and define $\Omega(0,\bar{q})$ to be the set of nspc.f.d. curves $\gamma\colon[0,\bar{x}_1]\longrightarrow\R^{2n+1}$ parameterized as in (\ref{AffCon}) and such that $\gamma(0)=0$, $\gamma(\bar{x}_1)=\bar{q}$. By Remark \ref{rem2}, for every $\gamma\in\Omega(0,\bar{q})$, $\pi\circ\gamma$ is contained in the set 
\begin{align*}
 \{(x_1,y_1,\dots,x_n,y_n,0):\; y_1^2+\sum_{j=2}^n(x_j^2+y_j^2)\leq x^2_1,\, 0 \leq x_1\leq \bar{x}_1\}\text{.}
\end{align*}
Using this and formulas (\ref{HorCon}), (\ref{AffCon}) we see that there exist constants $L>0$ and $M>0$ such that $|\dot{\gamma}(t)|\leq L$ and $|z(\gamma(t))|\leq M$ for every $\gamma\in\Omega(0,\bar{q})$, where $|\cdot|$ stands for the usual Euclidean norm in $\R^{2n+1}$, $z(\gamma)$ is the $z$-coordinate of $\gamma$, and $t\in[0,\bar{x}_1]$. 
Therefore the family of curves $\Omega(0,\bar{q})$ is equicontinuous and uniformly bounded. Since we are considering the trajectories of a control system whose set of control parameters is compact and convex, we also know that $\Omega(0,\bar{q})$ is closed with respect to uniform convergence (see \cite{Bress}). This establishes the following lemma.
\begin{lemma}\label{lem1}
 The set $\Omega(0,\bar{q})$ is compact with respect to the topology of uniform convergence. \qed
\end{lemma}
We now address the problem of the existence of maximizing geodesics. To this end, for a point $q\in M$, let us define the set 
\[
 \mathcal{D}^{tfd}_{q} = \{\lambda\in\mathcal{D}_q:\;\mathcal{H}(q,\lambda)<0,\,\langle\lambda,X_1(q)\rangle<0\} \text{.}
\]

\begin{proposition}
 For every $q_0\in\R^{2n+1}$ and $q\in J^+(q_0)$ there exists a maximizer joining $q_0$ to $q$. In the case where $q\in I^+(q_0)$, every such maximizer is a t.f.d. Hamiltonian geodesic. Consequently, 
 \begin{align}\label{exptim}
 \exp_{q_0}(\mathcal{D}^{tfd}_{q_0}) = I^+(q_0)\text{.}
 \end{align}
 On the other hand, when $q\in \partial J^+(q_0)$, every maximizer between $q_0$ and $q$ is a finite concatenation of null f.d. Hamiltonian geodesics.
\end{proposition}
\begin{proof}
Take two point $q_0,q\in M$ such that $q\in J^+(q_0)$. By applying an appropriate isometry we can assume that $q_0=0$ and $q\in J^+(0)$. By the above considerations we deduce that 
 \[
 L = \sup\{L(\gamma):\;\gamma\in\Omega(0,q)\} 
 \]
 is finite. That being so, we can find a sequence $\{\gamma_{\nu}\}\subset\Omega(0,q)$ such that $L(\gamma_{\nu})\longrightarrow L$. Using Lemma \ref{lem1} and \cite{Bress} or \cite{Gro2002} we make sure that the limit curve of $\{\gamma_{\nu}\}$ is nspc.f.d., and hence is a maximizing geodesic. In this way the existence of maximizing geodesics has been established.

Suppose now that $q\in I^+(q_0)$. By Section \ref{SecPrelimSL} we know that every maximizing geodesic connecting $q_0$ to $q$ is a t.f.d. Hamiltonian geodesic, thus proving \eqref{exptim}.

Finally, the last claim follows from \cite{Gro2009}.
\end{proof}

\subsection{Geometrically optimal curves are smooth}\label{subssmooth}
Recall that we are interested in the set $J^+(0)$ and its boundary. As it was pointed out earlier, every geometrically optimal curve is a concatenation of a finite number of null f.d. Hamiltonian geodesics. In this subsection we extract more detailed information about such curves. 

To begin with, consider the curve $\gamma:[0,+\infty)\longrightarrow\R^{2n+1}$ defined by
\[
 \gamma(t) = (a_1t,b_1t,\dots,a_nt,b_nt,0)\text{,}
\]
with $-a_1^2+b_1^2 + \sum_{j=2}^{n}(a_j^2+b_j^2)=0$, $a_1>0$. The curve $\gamma$ takes values in $\{z=0\}$ and, as is well known, is a maximizing null f.d. geodesic in the Minkowski space $(\R^{2n},g_0)$. Equivalently, in the terminology of Lorentzian geometry, $\gamma$ lies on the boundary of the chronological future of the origin in $(\R^{2n},g_0)$. 
Fix $t>0$ and suppose that $\eta:[0,T]\longrightarrow\R^{2n+1}$ is a null f.d. curve joining the origin to $\gamma(t)$. Then $\pi\circ\eta$ is a null f.d. curve in $(\R^{2n},g_0)$ joining $\gamma(0)$ to $\gamma(t)$. However, up to reparameterization, the only nspc.f.d. curve in $(\R^{2n},g_0)$ joining these two points is $\gamma_{|[0,t]}$. Hence, after a suitable reparameterization of $\eta$, $\gamma_{|[0,t]} = \pi\circ\eta$. The horizontality condition \eqref{HorCon} then implies $\gamma_{|[0,t]} = \eta$. Thus, for every $t>0$, $\gamma_{|[0,t]}$ is, up to reparameterization, the unique null f.d. curve connecting $\gamma(0)$ to $\gamma(t)$. Consequently, $\gamma$ is globally geometrically optimal.

Therefore, in this and the following subsections, we focus on geometrically optimal curves joining the origin to a point with nonzero $z$-coordinate.

Fix a point 
\[
 \bar{q}=(\bar{x}_1,\bar{y}_1,\dots,\bar{x}_n,\bar{y}_n,0)\in\pi(I^+(0))
\]
and consider the straight line 
\begin{align}\label{Linegamma}
 \Gamma = \{(\bar{x}_1,\bar{y}_1,\dots,\bar{x}_n,\bar{y}_n,z):\; z\in\R\} \text{.}
\end{align}
Denote by $\Omega(\Gamma)$ the set of nspc.f.d. curves $\gamma\colon[0,\bar{x}_1]\longrightarrow\R^{2n+1}$ such that $\gamma$ is parameterized as in (\ref{AffCon}), $\gamma(0)=0$, and $\gamma(\bar{x}_1)\in\Gamma$. Using reasoning similar to that which led us to Lemma \ref{lem1}, we can show that $\Omega(\Gamma)$ is compact with respect to topology of uniform convergence. 

Consider the mapping $F\colon\Omega(\Gamma)\longrightarrow\R$ defined by
\[
 F(\gamma) = z(\gamma(\bar{x}_1))
\]
or, equivalently, by
\[
 F(\gamma) = \tfrac12\int\displaylimits_{\tilde{\gamma}}\sum_{i=1}^n(y_idx_i - x_idy_i) \text{,}
\]
where $\tilde{\gamma} = \pi\circ\gamma$.
Evidently, $F$ is continuous and hence bounded, which ensures the existence of $\gamma_{\max},\gamma_{\min}\in\Omega(\Gamma)$ such that 
\begin{align}\label{zmax}
 z_{\max}=z(\gamma_{\max}(\bar{x}_1)) = \max_{\gamma\in\Omega(\Gamma)} F(\gamma)
\end{align}
and
\begin{align}\label{zmin}
 z_{\min}=z(\gamma_{\min}(\bar{x}_1)) = \min_{\gamma\in\Omega(\Gamma)} F(\gamma) \text{.}
\end{align}
Moreover, for a point $\bar{q}$ chosen as above, it is clear that
\[
 z_{\min}<0<z_{\max}\text{.}
\]
It is easy to see that $\gamma_{\max}(\bar{x}_1),\gamma_{\min}(\bar{x}_1)\in\partial J^+(0)$ and, consequently, $\gamma_{\max},\gamma_{\min}\subset\partial J^+(0)$. We can rephrase this by saying that $\gamma_{\max}$ and $\gamma_{\min}$ are geometrically optimal or, to put it another way, that $\gamma_{\max}$ (resp. $\gamma_{\min}$) is a maximizing geodesic joining $0$ to $\gamma_{\max}(\bar{x}_1)$ (resp. $\gamma_{\min}(\bar{x}_1)$). 

We will prove the stronger assertion that the curves obtained in this way are the only geometrically optimal curves starting from zero. Indeed, suppose that $z_1,z_2$ are two numbers satisfying $z_{\min}<z_1<z_2<z_{\max}$ and such that 
\[
 (\bar{x}_1,\bar{y}_1,\dots,\bar{x}_n,\bar{y}_n,z_i)\in\ I^+(0)\text{,}
\]
$i=1,2$. 
\begin{lemma}\label{Lem2}
 Under the above assumptions
 \[
  \{(\bar{x}_1,\bar{y}_1,\dots,\bar{x}_n,\bar{y}_n,z):\;z\in[z_1,z_2]\}\subset I^+(0)\text{.}
 \]
\end{lemma}
\begin{proof}
 Let $z\in(z_1,z_2)$ and choose t.f.d. curves $\gamma_i\in\Omega(\Gamma)$ such that $z_i = F(\gamma_i)$, $i=1,2$. Since $\pi(I^+(0))$ is convex, the curve
 \[
  \tilde{\gamma}_s = s \pi\circ\gamma_1 + (1-s)\pi\circ\gamma_2\text{,}
 \]
$s\in[0,1]$, lies in $\pi(I^+(0))$ and is t.f.d. in the sense of $(\R^{2n},g_0)$ -- cf. Remark \ref{rem2}. Denote by $\gamma_s$ the horizontal lift of $\tilde{\gamma}_s$ starting from zero. Evidently, $\gamma_s$ is t.f.d. and $\gamma_s\in\Omega(\Gamma)$ for every $s$. The function $[0,1] \ni s \longrightarrow F(\gamma_s)$ is continuous and attains the values $z_1$ and $z_2$, at $s=1$ and $s=0$, respectively. Thus, necessarily, $F(\gamma_s)=z$ for some $s\in(0,1)$. This proves that
\[
(\bar{x}_1,\bar{y}_1,\dots,\bar{x}_n,\bar{y}_n,z)\in I^+(0)\text{.}
\]
\end{proof}
Arguing as in the proof of Lemma \ref{Lem2}, but with $z_1=z_{\min}$ and $z_2=z_{\max}$, and remembering that for any null f.d. vectors $v_1,v_2$ and $s\in(0,1)$, the vector 
\[
 sv_1+(1-s)v_2
\]
is t.f.d. if and only if $v_1$ and $v_2$ are not proportional, we obtain the following result.
\begin{proposition}\label{propositionProj}
 \[
 \{(\bar{x}_1,\bar{y}_1,\dots,\bar{x}_n,\bar{y}_n,z):\; z\in [z_{\min},z_{\max}]\} \subset J^+(0)\text{,}
\]
and
 \[
 \{(\bar{x}_1,\bar{y}_1,\dots,\bar{x}_n,\bar{y}_n,z):\; z\in (z_{\min},z_{\max})\} \subset I^+(0)\text{.}
\]
Consequently, the projection $\pi$ induces bijections
\[
 \pi:\partial J^+(0)\cap\{z>0\} \longrightarrow \pi(I^+(0))\text{,}
\]
and
\[
 \pi:\partial J^+(0)\cap\{z<0\} \longrightarrow \pi(I^+(0))\text{.}
\]
\end{proposition}
\begin{proof}
Indeed, take $z\in (z_{\min},z_{\max})$. As established above, there exists $s\in (0,1)$ such that $F(\gamma_s)=z$, where $\gamma_s$ is the horizontal lift from the origin of the curve
\[
  \tilde{\gamma}_s = s \pi\circ\gamma_{\max} + (1-s)\pi\circ\gamma_{\min}\text{.}
\]
Suppose, to the contrary, that 
\[
 q=(\bar{x}_1,\bar{y}_1,\dots,\bar{x}_n,\bar{y}_n,z)\notin I^+(0)\text{.}
\]
Then $\gamma_s$ is null f.d. By the preceding observations, $\dot{\gamma}_{\max}(t)$ and $\dot{\gamma}_{\min}(t)$ must be proportional for almost every $t$. In view of (\ref{AffCon}), this is possible only if $\dot{\gamma}_{\max}(t) = \dot{\gamma}_{\min}(t)$ almost everywhere. Hence $\gamma_{\max} = \gamma_{\min}$, which is a contradiction.
\end{proof}

In order to describe the curves $\gamma_{\max},\gamma_{\min}$ we use PMP. We will solve the following optimization problem:\smallskip

\noindent \textit{Fix a point $(\bar{x}_1,\bar{y}_1,\dots,\bar{x}_n,\bar{y}_n,0)\in\pi(I^+(0))$ and the straight line \eqref{Linegamma}. Among all curves $\gamma\in\Omega(\Gamma)$ find ones for which the value $F(\gamma)$ is extremal.} \smallskip

If $\gamma:[0,\bar{x}_1]\longrightarrow \R^{2n+1}$ is parameterized as in \eqref{AffCon}, it is of the form 
\[
 \gamma(t)=(t,y_1(t),x_2(t),y_2(t),\dots,x_n(t),y_n(t),z(t))\text{,}
\]
and $F(\gamma)$ becomes
\begin{align*}
 F(\gamma) = \int_0^{\bar{x}_1}y_1dt - \frac12\bar{x}_1\bar{y}_1 + \frac12\int\displaylimits_{\pi\circ\gamma}\sum_{j=2}^n(y_jdx_j - x_jdy_j)\text{.}
\end{align*}
Therefore, we can reformulate our problem as follows. 

\begin{problem}\label{OptPr1}
 Let us fix a point $(\bar{y}_1,\bar{x}_2,\bar{y}_2,\dots,\bar{x}_n,\bar{y}_n)$ and a number $\bar{x}_1>0$ such that 
\[
 0<\bar{y}_1^2 + \sum_{j=2}^n(\bar{x}_j^2+\bar{y}_j^2) < \bar{x}_1^2.
\]
Among all curves $$[0,\bar{x}_1] \ni t \longrightarrow (y_1(t),x_2(t),y_2(t),\dots,x_n(t),y_n(t))$$ connecting the origin to $(\bar{y}_1,\bar{x}_2,\bar{y}_2,\dots,\bar{x}_n,\bar{y}_n)$, and such that
\[
 \dot{y}_1^2 + \sum_{j=2}^n(\dot{x}_j^2 + \dot{y}_j^2) = 1\text{,}
\]
find those for which the integral
\[
 \int\displaylimits_0^{\bar{x}_1}\big(y_1 + \frac12\sum_{j=2}^n(y_j\dot{x}_j - x_j\dot{y}_j)\big)dt 
\]
attains its maximal value.
\end{problem}

\begin{problem}\label{OptPr2}
 Let us fix a point $(\bar{y}_1,\bar{x}_2,\bar{y}_2,\dots,\bar{x}_n,\bar{y}_n)$ and a number $\bar{x}_1>0$ such that 
\[
 0<\bar{y}_1^2 + \sum_{j=2}^n(\bar{x}_j^2+\bar{y}_j^2) < \bar{x}_1^2.
\]
Among all curves $$[0,\bar{x}_1] \ni t \longrightarrow (y_1(t),x_2(t),y_2(t),\dots,x_n(t),y_n(t))$$ connecting the origin to $(\bar{y}_1,\bar{x}_2,\bar{y}_2,\dots,\bar{x}_n,\bar{y}_n)$, and such that
\[
 \dot{y}_1^2 + \sum_{j=2}^n(\dot{x}_j^2 + \dot{y}_j^2) = 1\text{,}
\]
find those for which the integral
\[
 \int\displaylimits_0^{\bar{x}_1}\big(y_1 + \frac12\sum_{j=2}^n(y_j\dot{x}_j - x_j\dot{y}_j)\big)dt 
\]
attains its minimal value.
\end{problem}

Remark that if $(y_1(t),x_2(t),y_2(t),\dots,x_n(t),y_n(t))$ is a solution to one of the above problems, then the corresponding geometrically optimal curve is given as the horizontal lift, starting from the origin, of the curve
\[
 \tilde{\gamma}(t)=(t,y_1(t),x_2(t),y_2(t),\dots,x_n(t),y_n(t)) \text{.}
\]

We solve Problem \ref{OptPr1}. In order to apply PMP, notice that the dynamics is given by
\begin{gather}
 \dot{y}_1= v_1,\;\dot{x}_j = u_j,\; \dot{y}_j = v_j,
\end{gather}
$j=2,\dots,n$, where the control parameters satisfy $v_1^2 + \sum_{j=2}^n(u_j^2 + v_j^2) = 1$, admissible controls are of class $L^{\infty}[0,\bar{x}_1]$, and the cost equals
\begin{gather}
 \int_0^{\bar{x}_1}\big(y_1 + \frac12\sum_{j=2}^n(y_ju_j - x_jv_j)\big)dt\text{.}
\end{gather}

Introduce the notation $w=(v_1,u_2,v_2,\dots,u_n,v_n)$, $\lambda = (q_1,p_2,q_2,\dots,p_n,q_n)$, $z=(y_1,x_2,y_2,\dots,x_n,y_n)$. The PMP Hamiltonian is
\begin{gather}
 h_w(z,\lambda) =
 q_1v_1 + \sum_{j=2}^n(p_ju_j+q_jv_j) - \nu\big(y_1 + \frac12\sum_{j=2}^n(y_ju_j - x_jv_j)\big) \text{.}
\end{gather}
We consider two cases.

\noindent \textbf{Abnormal case}. In this case we have $\nu=0$. The Hamiltonian is
\begin{gather*}
 h_w(z,\lambda) = q_1v_1 + \sum_{j=2}^n(p_ju_j+q_jv_j)\text{.}
\end{gather*}
Suppose that $(z(t),\lambda(t))$ is the corresponding bi-extremal, with optimal control $w(t)$. Then, equations $\dot{p}_j=-\frac{\partial h_{w(t)}}{\partial x_j}$, $j=2,\dots,n$, $\dot{q}_i=-\frac{\partial h_{w(t)}}{\partial y_i}$, $i=1,\dots,n$, reduce to
\begin{align}
 \dot{q}_1 = 0, \quad \dot{p}_j =0, \quad \dot{q}_j = 0 \text{,}
\end{align}
i.e., $q_1=q_{10},p_j=p_{j0},q_j=q_{j0}$ are constant, $j=2,\dots,n$. The PMP maximization condition, i.e.,
\[
 q_{10}v_1(t) + \sum_{j=2}^n(p_{j0}u_j(t)+q_{j0}v_j(t)) = \max_w\big(q_{10}v_1 + \sum_{j=2}^n(p_{j0}u_j+q_{j0}v_j)\big)\text{,}
\]
together with the Lagrange multipliers, shows that $v_1,u_j,v_j$, $j=2\dots$, are constant. It follows that the corresponding geometrically optimal curve would be contained in the plane $\{z=0\}$ which is impossible. Thus abnormal extrelmals of our problem do not exist.\smallskip

\noindent \textbf{Normal case}. In this case we have $\nu=-1$. The Hamiltonian is
\begin{gather*}
 h_w(z,\lambda) = q_1v_1 + \sum_{j=2}^n(p_ju_j+q_jv_j) + y_1 + \frac12\sum_{j=2}^n(y_ju_j - x_jv_j) \text{,}
\end{gather*}
and the corresponding Hamiltonian equations yield
\begin{gather}\label{solham}
 q_1 = q_{10}-t,\,\dot{p}_j = \frac12v_j, \, \dot{q}_j = -\frac12u_j, j=2,\dots,n \text{.}
\end{gather}

\textit{We adopt the following convention: throughout this and the following sections, we assume that $j$ ranges from $2$ to $n$.}

Let us rewrite the Hamiltonian $h_w$ in a more convenient form:
\begin{align*} 
 h_w(z,\lambda) = q_1v_1 + y_1 + \sum_j(P_ju_j+Q_jv_j) \text{,}
\end{align*}
where $P_j = p_j + \frac12y_j$, $Q_j = q_j - \frac12x_j$.
Suppose that $(z(t),\lambda(t))$ is a bi-extremal with optimal control $w(t)$. Applying the maximization condition of PMP
\begin{gather}\label{MaxPMP}
 h_{w(t)}(z(t),\lambda(t)) = \max_w h_w(z(t),\lambda(t)) \equiv const \text{,}
\end{gather}
and the Lagrange multipliers rule, we obtain
\begin{align}\label{LagMult}
 \begin{array}{l}
  q_1(t) = \mu(t) v_1(t) \smallskip \\
  P_j(t) = \mu(t)u_j(t)  \smallskip \\
  Q_j(t) = \mu(t)v_j(t)
 \end{array}\text{,}
\end{align}
from which
\begin{align}\label{mu}
 q_1(t)^2 + \sum_j(P_j(t)^2 + Q_j(t)^2) = \mu^2(t) \text{.}
\end{align}
Using \eqref{solham} and the definition of $P_j$ and $Q_j$ we have
\[
 \dot{P}_j = \dot{y}_j \quad\text{and}\quad \dot{Q}_j = -\dot{x}_j \text{.}
\]
Combining this with \eqref{LagMult} we find that 
\begin{align*}
 \frac{d}{dt}(P_j^2+Q_j^2) = 2\mu(t)(\dot{x}_j\dot{y}_j - \dot{y}_j\dot{x}_j)=0\text{.}
\end{align*}
It follows that $P_j(t)^2+Q_j(t)^2$ is constant. Suppose, for a moment, that 
\begin{align}\label{PQzero}
 \sum_j(P_j(t)^2 + Q_j(t)^2) = 0\text{.}
\end{align}
Then $P_j(t)=Q_j(t)=0$ identically for every $j$. Moreover, by \eqref{solham}
\[
 \dot{p}_j= \frac12\dot{y}_j \text{,}
\]
whereas the definition of \(P_j\) gives
\[
 \dot{p}_j= -\frac12\dot{y}_j \text{.}
\]
Therefore, $\dot{y}_j=0$, and hence $y_j(t) = 0$ identically for every $j$. A similar argument shows that $x_j(t)= 0$ identically for every $j$. It follows that the resulting geometrically optimal curve is contained in $\{z=0\}$, a possibility that has already been excluded. Thus \eqref{PQzero} cannot hold and, consequently, $\mu(t)$ never vanishes.

Now, using the observations above, we obtain the system
\begin{align}\label{EqPQ}
 \begin{aligned}
  \dot{P}_j &= \dfrac{Q_j}{\mu} \\
  \dot{Q}_j &= -\dfrac{P_j}{\mu}
 \end{aligned} \text{.}
\end{align}
Setting $R_j = \sqrt{P_j(t)^2+Q_j(t)^2}$, we can write $P_j(t) = R_j\cos\alpha_j(t)$, $Q_j(t) = R_j\sin\alpha_j(t)$. Substituting it into \eqref{EqPQ} we obtain
\begin{align*} 
 \begin{aligned} 
   P_j(t) &= R_j\cos\left( \varphi_j-\int_0^t\frac{d\tau}{\mu(\tau)} \right),\\ 
   Q_j(t) &= R_j\sin\left( \varphi_j-\int_0^t\frac{d\tau}{\mu(\tau)} \right), 
 \end{aligned} 
\end{align*}
for suitably chosen constants $\varphi_j$, which finally leads to 
\begin{align}\label{XJYJ} 
 \begin{aligned} 
  x_j(t) &= -R_j\left( \sin\left( \varphi_j-\int_0^t\frac{d\tau}{\mu(\tau)} \right) -\sin\varphi_j \right),\\ 
  y_j(t) &= R_j\left( \cos\left( \varphi_j-\int_0^t\frac{d\tau}{\mu(\tau)} \right) -\cos\varphi_j \right). 
 \end{aligned} 
\end{align}

Next, directly from \eqref{mu} we have
\begin{align}\label{MU}
 \mu(t) = \pm\sqrt{(q_{10}-t)^2 + \sum\nolimits_jR_j^2}\text{.}
\end{align}
Substituting \eqref{LagMult} and \eqref{MU} into the maximization condition \eqref{MaxPMP} we obtain
\[
 h_{w(t)}(z(t),\lambda(t)) = y_1(t) + \mu(t) = \mu(0)\text{.}
\]
It follows that $\mu$ is positive, i.e.,
\begin{align*}
 \mu(t) = \sqrt{(q_{10}-t)^2 + \sum_jR_j^2} \text{.}
\end{align*}
The preceding identity also yields the equation for $y_1$. Thus, the desired geometrically optimal curve corresponding to the extremal just obtained is given by the horizontal lift, starting from the origin, of the curve defined by \eqref{XJYJ} together with
\begin{align}\label{X1Y1}
 \begin{array}{l}
  x_1(t) = t \bigskip \\
  y_1(t) = -(\mu(t) - \mu(0))
 \end{array}\text{.}
\end{align}
In particular, this curve is smooth.

To solve Problem \ref{OptPr2} we proceed similarly, except that in this case $\mu$ is of the form
\begin{align*}
 \mu(t) = -\sqrt{(q_{10}+t)^2 + \sum_jR_j^2}
\end{align*}
and $\nu = 1$.

\subsection{Geometrically optimal curves are null f.d. Hamiltonian geodesics}
Since it is not hard to construct, either in sub-Lorentzian or in sub-Riemannian geometry, a smooth curve that is a concatenation of distinct Hamiltonian geodesics, the claim stated in the title of this subsection requires proof.

To this end, let us recall equations \eqref{GeodEqsrnotzero}. As we are interested in null f.d. geodesics starting from zero and not contained in $\{z=0\}$, we take an initial covector of the form
\begin{align}\label{incon}
 A_1 = -1, B_1 = \pm\sqrt{1-\sum_j(A_j^2+B_j^2)},A_2,B_2,\dots,A_n,B_n,r \text{,}
\end{align}
where $|B_1|<1$ and $r \ne 0$.
Define $S_j = \sqrt{A_j^2+B_j^2}$ and write 
\begin{align}\label{SJ}
 A_j = (\mathrm{sgn\,}r)\,S_j\cos\theta_j, \quad B_j = (\mathrm{sgn\,}r)\,S_j\sin\theta_j\text{.}
\end{align}
Substituting \eqref{incon} and \eqref{SJ} to \eqref{GeodEqsrnotzero} we obtain
\begin{align}\label{NullGeod}
\begin{array}{l}
 x_1(t) = \frac{1}{r}\sinh(rt) - \frac{B_1}{r}(\cosh(rt)- 1) \bigskip \\
 y_1(t) = \frac{B_1}{r}\sinh(rt) - \frac{1}{r}(\cosh(rt) -1)  \bigskip \\  
 x_j(t) = -\frac{S_j}{|r|}\Big(\sin(\theta_j-rt) - \sin\theta_j\Big) \bigskip \\
 y_j(t) = \frac{S_j}{|r|}\Big(\cos(\theta_j-rt) - \cos\theta_j\Big)
\end{array}
\end{align}
Solving the first equation in \eqref{NullGeod} with respect to $t$ we get
\begin{align}\label{x1solved}
 t = t(s) = \frac{1}{r}\ln\frac{B_1-rs-\sqrt{r^2s^2 - 2B_1rs+1}}{B_1-1} \text{,}
\end{align}
where we write $s$ instead of $x_1$.

Now, we substitute \eqref{x1solved} to \eqref{NullGeod}. 
By computing $\frac{d}{ds}(y_1(t(s)))$, integrating the resulting expression, and performing some algebraic manipulations, we see that the first two equations in \eqref{NullGeod} take the following form:
\begin{align}\label{FirstTwo}
\begin{array}{l}
 x_1(s) = s \bigskip \\
 y_1(s) = -(\mu(s) - \mu(0)) \text{,}
\end{array}
\end{align}
where
\begin{align*}
 \mu(s) = \pm\sqrt{(\tfrac{B_1}{r}-s)^2 + \tfrac{1-B_1^2}{r^2}}\text{.}
\end{align*}
The plus sign corresponds to $r>0$, whereas the minus sign corresponds to $r<0$.
Here and below by $x_1(s)$ we mean $x_1(t(s))$, and similarly for the remaining coordinates.
Finally, let us look at the equations for $x_j$ and $y_j$. By direct calculations we make sure that
\begin{align}
\int\displaylimits_0^s\frac{d\tau}{\mu(\tau)} = rt(s) \text{.}
\end{align}
Now we conclude that equations \eqref{NullGeod}, after replacing $t$ with $t(s)$, coincide with the extremals either of Problem \ref{OptPr1} or Problem \ref{OptPr2}, where 
\[
 R_j = \tfrac{S_j}{|r|},\quad \varphi_j=\theta_j,\quad q_{10} = \pm\frac{B_1}{r}\text{.}
\]

In this we we have established Theorem \ref{Th2}.

\subsection{The boundary of the reachable set}
By Proposition \ref{propositionProj}, the set $J^+(0)$ is convex in the $z$-direction. This allows us to state the following result.
\begin{proposition}\label{proBound}
 There exist functions 
 \[
  \varphi^+,\varphi^-\colon \pi(J^+(0))\longrightarrow\R
 \]
such that
 \begin{align*}
  J^+(0)\cap\{z\geq 0\} = \{q\in\R^{2n+1}:\; \pi(q)\in\pi(J^+(0)),\; 0\leq z \leq \varphi^+(\pi(q))\}
 \end{align*}
and
\begin{align*}
  J^+(0)\cap\{z\leq 0\} = \{q\in\R^{2n+1}:\; \pi(q)\in\pi(J^+(0)),\; \varphi^-(\pi(q))\leq z \leq 0\}\text{.}
 \end{align*}
 Consequently, 
 \begin{align*}
  \partial J^+(0)\cap\{z>0\} = \{q\in\R^{2n+1}:\; \pi(q)\in\pi(I^+(0)),\; z = \varphi^+(\pi(q))\}
 \end{align*}
and
\begin{align*}
  \partial J^+(0)\cap\{z<0\} = \{q\in\R^{2n+1}:\; \pi(q)\in\pi(I^+(0)),\; z = \varphi^-(\pi(q))\}\text{,}
 \end{align*}
 where $q=(x_1,y_1,\dots,x_n,y_n,z)$. \qed
\end{proposition}

We will need the following well-known fact from linear Lorentzian geometry.
\begin{lemma}\label{NSPCVEC}
 Let $V$ be a vector space equipped with the scalar product $g$ of Lorentzian signature. We assume that $V$ is time-oriented. Suppose that $0\ne v\in V$ is such that $g(v,w)\leq 0$ for every nspc.f.d. vector $w\in V$. Then $v$ is nspc.f.d. If, in addition, there exists a null f.d. vector $w\in V$ such that $g(v,w)=0$, then $v$ is also null f.d. and is proportional to $w$.
\end{lemma}
\begin{proof}
 Suppose the contrary, i.e., that $v$ is spacelike. Then the orthogonal complement $v^{\perp}$ has Lorentzian signature. Choose a t.f.d. $u\in v^{\perp}$. For a sufficiently small $\varepsilon>0$, the vector $u+\varepsilon v$ is t.f.d. Therefore, by assumption, $g(v,u+\varepsilon v)\leq 0$. However, by the choice of $u$,
 \[  
   g(v,u+\varepsilon v) = \varepsilon g(v,v)>0
 \]
which is a contradiction. It follows that $v$ is nspc.f.d., and the second claim is now clear.
\end{proof}

Introduce the following functions:
\[
\Phi^{\pm}(x_1,y_1,\dots,x_n,y_n,z) = z - \varphi^{\pm}(x_1,y_1,\dots,x_n,y_n) \text{.}
\]
The preceding discussion yields the following proposition. 
\begin{proposition}\label{PropBound}
 Suppose that $\varphi^{+}$ is smooth on some open set $U$ in $\R^{2n}$, where as above we identify $\R^{2n}$ with the subspace $\{z=0\}$ in $\R^{2n+1}$. Then $\nabla_H\Phi^{+}$ is null f.d. when restricted to the set $\partial J^+(0)\cap \pi^{-1}(U)$. Moreover, if $\gamma$ is a geometrically optimal curve that starts from the origin and enters $\pi^{-1}(U)$, then $\gamma\cap \pi^{-1}(U)$ is a trajectory of $\nabla_H\Phi^{+}$. Consequently, if $[a,b]$ is an interval such that $\gamma([a,b])\subset \pi^{-1}(U)$, then $\gamma_{|[a,b]}$ is the unique maximizer contained in $\pi^{-1}(U)$ that joins $\gamma(a)$ to $\gamma(b)$. The analogous statement holds for $\varphi^-$ and $\Phi^-$.
\end{proposition}
\begin{proof}
 Suppose that $\varphi^{+}$ is smooth on some open set $U\subset\R^{2n}$. Then $\Phi^{+}$ is smooth on $\pi^{-1}(U)$ and clearly $\nabla_H\Phi^+\ne 0$. Moreover, by Proposition \ref{proBound},
\[
  \partial J^+(0)\cap\pi^{-1}(U) = \{\Phi^+=0\}\cap\pi^{-1}(U)\text{.}
\]
Denote this set by $N$. It is a smooth hypersurface. Take an arbitrary point $q\in N$ and a nspc.f.d. $z\in T_qN\cap H_q$. Extend $z$ to a smooth nspc.f.d. vector field $Z$ defined in a neighborhood of $q$, and consider the trajectory $
\sigma\colon[0,T]\longrightarrow\pi^{-1}(U)$ of $Z$ with $\sigma(0)=q$. Since $\sigma([0,T])\subset J^+(0)$, we must have 
\begin{align*}
 \partial_Z(\Phi^+) = g(Z,\nabla_H\Phi^+)\leq 0
\end{align*}
at $q$, so in particular 
\begin{align*}
 g(z,\nabla_H\Phi^+(q))\leq 0 \text{.}
\end{align*}
Furthermore, $N$ is formed by geometrically optimal curves which, as we already know, are null f.d. and smooth. Let $\gamma:[a,b]\longrightarrow N$ be one such curve passing through $q$, so that $\gamma(t_0)=q$ for some $t_0\in(a,b)$. Since $\Phi^+(\gamma(t))=0$ identically,
\[
 g(\dot{\gamma}(t_0),\nabla_H\Phi^+(\gamma(t_0))) = 0 \text{.}
\]
By Lemma \ref{NSPCVEC}, this shows that $\nabla_H\Phi^+(q)$ is a null f.d. vector proportional to $\dot{\gamma}(t_0)$. Since $q$ was arbitrary, it follows that the restriction of $\nabla_H\Phi^+$ to $\partial J^+(0)\cap\pi^{-1}(U)$ is null f.d. Moreover, $\nabla_H\Phi^+(\gamma(t))$ is proportional to $\dot{\gamma}(t)$ for every $t$. Thus, up to a change of parameter, $\gamma$ is a trajectory of $\nabla_H\Phi^+$. The uniqueness of $\gamma$ now follows from \cite{Gro2009}.
\end{proof}

\subsection{Complete characterization of geometrically optimal curves}\label{CompCharGeomOpt}
Let 
\begin{align}\label{TheSetD}
\begin{array}{l}
\mathcal{D} = \{(-1,B_1,A_2,B_2,\dots,A_n,B_n,r,t) :\smallskip  \\
 B_1^2 + \sum\nolimits_j(A_j^2+B_j^2) = 1,\, -1<B_1<1,\, r\ne 0,\,t\in (0,\tfrac{2\pi}{|r|}) \}\subset T_0^*\R^{2n+1}\times(0,\tfrac{2\pi}{|r|}) 
\end{array}\text{.}
\end{align}
Topologically, $\mathcal{D}$ is
\[
 S^{2n-2}\backslash\{p_{-1},p_{+1}\}\times\R\backslash\{0\}\times (0,\tfrac{2\pi}{|r|}) \text{,}
\]
where $S^{2n-2}$ is the unit $(2n-2)$-dimensional sphere, and $p_{\pm 1}$ is a point of $S^{2n-2}$ corresponding to $B_1=\pm 1$.
Consider the mapping 
\begin{align}\label{Map}
 \mathcal{D} \ni (-1,B_1,A_2,B_2,\dots,A_n,B_n,r,t)\longrightarrow (x_1(t),y_1(t),\dots,x_n(t),y_n(t),z(t)) \in\R^{2n+1} \text{,}
\end{align}
where $x_i(t)$, $y_i(t)$, and $z(t)$ are given by \eqref{GeodEqsrnotzero}. As above we suppress the dependence on the initial covector and write $x_1(t)$, etc., instead of $x_1(t;-1,B_1,A_2,B_2,\allowbreak\dots,A_n,B_n,r)$, etc.

Our first aim in this subsection is to show that the image of the map \eqref{Map} is a smooth (in fact, analytic) hypersurface contained in $\partial J^+(0)\cap\{z\ne 0\}$. Note that if we take $B_1=\pm 1$, then the corresponding geodesic is a half-line issuing from zero and contained in $\{z=0\}$, and that is why we exclude this case -- cf. the beginning of Subsection \ref{subssmooth}.

We will compute the rank of the map
\[
  F\colon\mathcal{D}\ni(-1,B_1,A_2,B_2,\dots,A_n,B_n,r,t)\longrightarrow (x_1(t),y_1(t),\dots,x_n(t),y_n(t))\in\R^{2n}\text{.}
\]
Therefore, consider the family of null f.d. geodesics starting from the origin, and parameterized by initial covector as in \eqref{TheSetD}:
\begin{align}\label{geodrow}
\left\{\begin{array}{l}
 x_1(t) = \frac{1}{r}\sinh(rt) - \frac{B_1}{r}(\cosh(rt)- 1) \bigskip \\
 y_1(t) = \frac{B_1}{r}\sinh(rt) - \frac{1}{r}(\cosh(rt) -1)  \bigskip \\  
 x_j(t) = \frac{A_j}{r}\sin(rt) - \frac{B_j}{r}(\cos(rt) -1) \bigskip \\
 y_j(t) = \frac{B_j}{r}\sin(rt) + \frac{A_j}{r}(\cos(rt) -1) \bigskip \\
 z(t) = \sum_j\frac{A_j^2+B_j^2}{2r^2}\sinh(rt) - \sum_j \frac{A_j^2+B_j^2}{2r^2}\sin(rt)
\end{array}\right.\text{.}
\end{align}

We begin with the case $B_1\ne0$. To this end we assume that $B_1 = \pm\sqrt{1-\sum_j(A_j^2+B_j^2)}$.
We will carry out all computations for 
\begin{align}\label{B1subs}
 B_1 = \sqrt{1-\sum\nolimits_j(A_j^2+B_j^2)}\text{,}
\end{align}
the case with the negative sign being analogous.
Denote by $\tilde{\mathcal{D}}$ the set 
\[
 \{(A_2,B_2,\dots,A_n,B_n,r,t)\in\R^{2n}:\; 0<\sum_j(A_j^2+B_j^2)<1,\,r\ne0,\, t\in(0,\tfrac{2\pi}{|r|})\}
\]
and define the mapping $\tilde{F}$ by
\[
 \tilde{F}\colon\tilde{\mathcal{D}}\ni(A_2,B_2,\dots,A_n,B_n,r,t)\longrightarrow (x_1(t),y_1(t),\dots,x_n(t),y_n(t))\in\R^{2n}\text{.}
\]
We now compute the Jacobian of $\tilde{F}$. Before doing so, let us explain the restriction $t\in(0,\tfrac{2\pi}{|r|})$ in the definition of $\tilde{\mathcal{D}}$. Looking at \eqref{geodrow} it is easy to see that any two geodesics $\gamma$, $\tilde{\gamma}$ starting from zero, with respective initial covectors 
\begin{align}\label{InCov1}
 (-1,\pm\sqrt{1-\sum\nolimits_j(A_j^2+B_j^2)},A_2,B_2,\dots,A_n,B_n,r)
\end{align}
and
\begin{align}\label{InCov2}
 (-1,\pm\sqrt{1-\sum\nolimits_j(\tilde{A}_j^2+\tilde{B}_j^2)},\tilde{A}_2,\tilde{B}_2,\dots,\tilde{A}_n,\tilde{B}_n,r)\text{,}
\end{align}
satisfy $\gamma(\frac{2\pi}{|r|}) = \tilde{\gamma}(\frac{2\pi}{|r|})$, provided that $\sum_j(A_j^2+B_j^2) = \sum_j(\tilde{A}_j^2+\tilde{B}_j^2)$. Consequently, the Jacobian of $\tilde{F}$ becomes singular at $t = \frac{2\pi}{|r|}$ which justifies excluding this value. 

We make one further convention. Here and throughout what follows, we perform the computations for $r>0$, that is, for geodesics that start from the origin and subsequently remain in the half-space $\{z>0\}$. The case $r<0$, corresponding to geodesics that enter the lower half-space $\{z<0\}$, is analogous.

In order to compute the Jacobian of $\tilde{F}$, it is enough to find the determinant of the matrix 
\begin{align*}
\tilde{M}= \begin{pmatrix}
  \frac{\partial x_2}{\partial A_2} & \frac{\partial y_2}{\partial A_2} & \frac{\partial x_3}{\partial A_2} & \frac{\partial y_3}{\partial A_2} & \hdots & \frac{\partial x_n}{\partial A_2} & \frac{\partial y_n}{\partial A_2} & \frac{\partial x_1}{\partial A_2} & \frac{\partial y_1}{\partial A_2} \smallskip \\
  \frac{\partial x_2}{\partial B_2} & \frac{\partial y_2}{\partial B_2} & \frac{\partial x_3}{\partial B_2} & \frac{\partial y_3}{\partial B_2} & \hdots & \frac{\partial x_n}{\partial B_2} & \frac{\partial y_n}{\partial B_2} & \frac{\partial x_1}{\partial B_2} & \frac{\partial y_1}{\partial B_2} \smallskip \\
  \frac{\partial x_2}{\partial A_3} & \frac{\partial y_2}{\partial A_3} & \frac{\partial x_3}{\partial A_3} & \frac{\partial y_3}{\partial A_3} &\hdots & \frac{\partial x_n}{\partial A_3} & \frac{\partial y_n}{\partial A_3} & \frac{\partial x_1}{\partial A_3} & \frac{\partial y_1}{\partial A_3} \smallskip \\
  \frac{\partial x_2}{\partial B_3} & \frac{\partial y_2}{\partial B_3} & \frac{\partial x_3}{\partial B_3} & \frac{\partial y_3}{\partial B_3} & \hdots & \frac{\partial x_n}{\partial B_3} & \frac{\partial y_n}{\partial B_3} & \frac{\partial x_1}{\partial B_3} & \frac{\partial y_1}{\partial B_3} \smallskip \\
  \vdots & \vdots & \vdots & \vdots & \ddots & \vdots & \vdots & \vdots & \vdots \smallskip \\
  \frac{\partial x_2}{\partial A_n} & \frac{\partial y_2}{\partial A_n} & \frac{\partial x_3}{\partial A_n} & \frac{\partial y_3}{\partial A_n} & \hdots & \frac{\partial x_n}{\partial A_n} & \frac{\partial y_n}{\partial A_n} & \frac{\partial x_1}{\partial A_n} & \frac{\partial y_1}{\partial A_n} \smallskip \\
  \frac{\partial x_2}{\partial B_n} & \frac{\partial y_2}{\partial B_n} & \frac{\partial x_3}{\partial B_n} & \frac{\partial y_3}{\partial B_n} & \hdots & \frac{\partial x_n}{\partial B_n} & \frac{\partial y_n}{\partial B_n} & \frac{\partial x_1}{\partial B_n} & \frac{\partial y_1}{\partial B_n} \smallskip \\
  \frac{\partial x_2}{\partial r} & \frac{\partial y_2}{\partial r} & \frac{\partial x_3}{\partial r} & \frac{\partial y_3}{\partial r} & \hdots & \frac{\partial x_n}{\partial r} & \frac{\partial y_n}{\partial r} & \frac{\partial x_1}{\partial r} & \frac{\partial y_1}{\partial r} \smallskip \\
  \frac{\partial x_2}{\partial t} & \frac{\partial y_2}{\partial t} & \frac{\partial x_3}{\partial t} & \frac{\partial y_3}{\partial t} & \hdots & \frac{\partial x_n}{\partial t} & \frac{\partial y_n}{\partial t} & \frac{\partial x_1}{\partial t} & \frac{\partial y_1}{\partial t}
 \end{pmatrix}\text{.}
\end{align*}\smallskip
Let us rewrite $\tilde{M}$ as a block matrix, with blocks of order $2\times 2$: 
\begin{align*}
 \tilde{M}=\begin{pmatrix}
  P_2 & \text{\large 0} & \text{\large 0} &\hdots  & w_2\cdot v^T\\
  \text{\large 0} & P_3 & \text{\large 0} & \hdots & w_3\cdot v^T\\
  \vdots & \vdots & \ddots & \vdots & \vdots\\
  \text{\large 0} & \text{\large 0} & \hdots & P_n & w_n\cdot v^T \\
  C_2 & \hdots & \hdots & C_n & C_1
   \end{pmatrix}\text{,}
\end{align*}
where
\begin{gather}\label{MatrPj}
P_j = \begin{pmatrix}
        \frac{\partial x_j}{\partial A_j} & \frac{\partial y_j}{\partial A_j} \smallskip \\
        \frac{\partial x_j}{\partial B_j} & \frac{\partial y_j}{\partial B_j}
      \end{pmatrix} = 
      \begin{pmatrix}
       \frac{\sin(rt)}{r} & \frac{\cos(rt)-1}{r} \smallskip \\
       -\frac{\cos(rt)-1}{r} & \frac{\sin(rt)}{r} 
      \end{pmatrix}\text{,}
\end{gather}
\begin{gather}\label{MatrC1}
 C_1 = \begin{pmatrix}
        \frac{\partial x_1}{\partial r} & \frac{\partial y_1}{\partial r} \smallskip  \\
        \frac{\partial x_1}{\partial t} & \frac{\partial y_1}{\partial t}
       \end{pmatrix} = \\       
       \begin{pmatrix}
        \tfrac{(rt + B_1)\cosh(rt) - (B_1rt + 1)\sinh(rt)  - B_1}{r^2} & \tfrac{(B_1rt + 1)\cosh(rt) - (B_1 + rt)\sinh(rt) - 1}{r^2}\smallskip \\
        \cosh(rt) - B_1\sinh(rt) & B_1\cosh(rt) - \sinh(rt)
       \end{pmatrix}\notag
       \text{,}
\end{gather}\smallskip
with $B_1$ as in \eqref{B1subs},
\begin{gather}\label{MatrCj}
 C_j = \begin{pmatrix}
        \frac{\partial x_j}{\partial r} & \frac{\partial y_j}{\partial r} \smallskip \\
        \frac{\partial x_j}{\partial t} & \frac{\partial y_j}{\partial t}
       \end{pmatrix} = 
       \begin{pmatrix}
        \tfrac{(A_jrt + B_j)\cos(rt) + (B_jrt - A_j)\sin(rt) - B_j}{r^2} & \frac{(B_jrt  - A_j)\cos(rt) - (A_jrt + B_j)\sin(rt) + A_j}{r^2} \\
        A_j\cos(rt) + B_j\sin(rt) & B_j\cos(rt) - A_j\sin(rt)
       \end{pmatrix}\text{,}
\end{gather}
and finally
\begin{align*}
 w_j= \begin{pmatrix}
       \frac{\partial B_1}{\partial A_j} \smallskip\\
       \frac{\partial B_1}{\partial B_j}
      \end{pmatrix} = 
      \begin{pmatrix}
       -\frac{A_j}{B_1} \smallskip\\
       -\frac{B_j}{B_1}
      \end{pmatrix}\text{,}
\end{align*}
\begin{align*}
 v= \begin{pmatrix}
     \frac{\partial x_1}{\partial B_1} \smallskip\\
     \frac{\partial y_1}{\partial B_1}
    \end{pmatrix} =
    \begin{pmatrix}
     -\frac{\cosh(rt)-1}{r} \smallskip \\
     \frac{\sinh(rt)}{r}
    \end{pmatrix}\text{,}
\end{align*}
and
\begin{align*}
 {\large 0}= \begin{pmatrix}
     0 & 0 \smallskip\\
     0 & 0
    \end{pmatrix}\text{.}
\end{align*}

Next we rewrite $\tilde{M}$ in yet another block form, namely
\begin{align}
 \tilde{M} = \begin{pmatrix}
      P & R \\
      S & C_1
     \end{pmatrix}\text{,}
\end{align}
where 
\begin{align*}
P=\begin{pmatrix}
  P_2 & \text{\large 0} & \hdots &\text{\large 0}\\
  \text{\large 0} & P_3 & \hdots & \text{\large 0}\\
  \vdots & \vdots & \ddots & \vdots & \\
  \text{\large 0} & \text{\large 0} & \hdots & P_n
   \end{pmatrix}
\end{align*}
is of order $2n-2\times 2n-2$,
\begin{align*}
 R = \begin{pmatrix}
      w_2\cdot v^T \\
      \vdots \\
      w_n\cdot v^T
     \end{pmatrix}
\end{align*}
is of order $2n-2\times 2$,
\begin{align*}
 S = \begin{pmatrix}
      C_2 & \hdots & C_n
     \end{pmatrix}
\end{align*}
is of order $2\times 2n-2$, and $C_1$ is of order $2\times 2$. Since the matrix $P$ is invertible, we can apply Schur's formula (\cite{Gant}) to compute the determinant:
\begin{align}
 \det \tilde{M} = \det P \det(C_1 - SP^{-1}R)\text{,}
\end{align}
where 
\begin{align}\label{DetPj}
 \det P = (\det P_2)^{n-1} = (-2)^{n-1}\frac{(\cos(rt)-1)^{n-1}}{r^{2n-2}}\text{.}
\end{align}
Further, notice that
\[
 C_1 - SP^{-1}R = C_1 - (\sum_jC_jP_j^{-1}w_j)\cdot v^T\text{,}
\]
and since 
\[
 \det C_1 = \frac{\cosh(rt)-1}{r^2} \sum_j(A_j^2+B_j^2) \ne 0 \text{,}
\]
we can use, for instance, rank-one updates idea (see e.g. \cite[Chapter 6]{Meyer}):
\begin{align*}
 \det(C_1 - SP^{-1}R) = \det C_1\cdot\big(1 - v^TC_1^{-1}\sum_jC_jP_j^{-1}w_j\big) \text{.}
\end{align*}
After carrying out the necessary computations, we find that the last expression equals
\[
 \frac{\sum_j(A_j^2+B_j^2)\big(\sinh(rt)(\cos(rt)-1)+\sin(rt)(\cosh(rt)-1)\big)}{r^2(\cos(rt)-1)\sqrt{1-\sum_j(A_j^2+B_j^2)}} \text{.}
\] 
Combining the preceding results yields the following formula for $\det\tilde{M}$:
\begin{align*}
 \frac{(-2)^{n-1}\sum_j(A_j^2+B_j^2)\big(\sinh(rt)(\cos(rt)-1)+\sin(rt)(\cosh(rt)-1)\big)(\cos(rt)-1)^{n-2}}{r^{2n}\sqrt{1-\sum_j(A_j^2+B_j^2)}}\text{.}
\end{align*}
Since 
\[
 \sinh(rt)(\cos(rt)-1)+\sin(rt)(\cosh(rt)-1)
\]
is nonzero for $r\ne 0$ and $t\in(0,\frac{2\pi}{|r|})$, it follows that the Jacobian of $\tilde{F}$ does not vanish on $\tilde{\mathcal{D}}$. \smallskip

We now consider the remaining case. Namely, we take the initial covectors as in \eqref{TheSetD}, allowing $B_1=0$. Then at least one of the numbers $A_2,B_2,\dots,A_n,B_n$ must be nonzero. Suppose, for instance, that $A_2\ne 0$; the cases $B_2\ne0$, $A_3\ne0$, and so on are treated analogously. We may then assume that
\begin{align}\label{A2}
 A_2 = \pm\sqrt{1 - B_1^2 - B_2^2 - \sum\nolimits_k(A_k^2+B_k^2)}\text{,}
\end{align}
where we agree that the index $k$ ranges from $3$ to $n$. Consider the set 
\begin{align*}
 \begin{aligned}
  &\tilde{\tilde{D}} = \Big\{(B_1,B_2,A_3,B_3\dots,A_n,B_n,r,t)\in\R^{2n}:\; \\
  &B_1^2+B_2^2+\sum\nolimits_k(A_k^2+B_k^2)<1,\,|B_1|<1,\,r\ne0,\, t\in(0,\tfrac{2\pi}{|r|})\Big\}
 \end{aligned}
\end{align*}
and the mapping
\begin{align}
\tilde{\tilde{F}}\colon\tilde{\tilde{\mathcal{D}}}\ni(B_1,B_2,A_3,B_3,\dots,A_n,B_n,r,t)\longrightarrow (x_1(t),y_1(t),\dots,x_n(t),y_n(t))\in\R^{2n}\text{,}
\end{align}
where, as before, $x_1(t)$, etc., are given by \eqref{geodrow}, and we additionally substitute \eqref{A2} for $A_2$. In what follows, we set 
\begin{align}\label{A2subs}
 A_2 = \sqrt{1 - B_1^2 - B_2^2 - \sum_k(A_k^2+B_k^2)}\text{.}
\end{align}

To compute the Jacobian of $\tilde{\tilde{F}}$, let us consider the matrix
\begin{align*}
\tilde{\tilde{M}}= \begin{pmatrix}
 \frac{\partial x_1}{\partial r} & \frac{\partial y_1}{\partial r} & \frac{\partial x_3}{\partial r} & \frac{\partial y_3}{\partial r} & \hdots & \frac{\partial x_n}{\partial r} & \frac{\partial y_n}{\partial r} & \frac{\partial x_2}{\partial r} & \frac{\partial y_2}{\partial r} \smallskip \\
  \frac{\partial x_1}{\partial t} & \frac{\partial y_1}{\partial t} & \frac{\partial x_3}{\partial t} & \frac{\partial y_3}{\partial t} & \hdots & \frac{\partial x_n}{\partial t} & \frac{\partial y_n}{\partial t} & \frac{\partial x_2}{\partial t} & \frac{\partial y_2}{\partial t} \smallskip \\
  \frac{\partial x_1}{\partial A_3} & \frac{\partial y_1}{\partial A_3} & \frac{\partial x_3}{\partial A_3} & \frac{\partial y_3}{\partial A_3} & \hdots & \frac{\partial x_n}{\partial A_3} & \frac{\partial y_n}{\partial A_3} & \frac{\partial x_2}{\partial A_3} & \frac{\partial y_2}{\partial A_3} \smallskip \\
  \frac{\partial x_1}{\partial B_3} & \frac{\partial y_1}{\partial B_3} & \frac{\partial x_3}{\partial B_3} & \frac{\partial y_3}{\partial B_3} & \hdots & \frac{\partial x_n}{\partial B_3} & \frac{\partial y_n}{\partial B_3} & \frac{\partial x_2}{\partial B_3} & \frac{\partial y_2}{\partial B_3} \smallskip \\
  \vdots & \vdots & \vdots & \vdots & \ddots & \vdots & \vdots & \vdots & \vdots \smallskip \\
  \frac{\partial x_1}{\partial A_n} & \frac{\partial y_1}{\partial A_n} & \frac{\partial x_3}{\partial A_n} & \frac{\partial y_3}{\partial A_n} & \hdots & \frac{\partial x_n}{\partial A_n} & \frac{\partial y_n}{\partial A_n} & \frac{\partial x_2}{\partial A_n} & \frac{\partial y_2}{\partial A_n} \smallskip \\
  \frac{\partial x_1}{\partial B_n} & \frac{\partial y_1}{\partial B_n} & \frac{\partial x_3}{\partial B_n} & \frac{\partial y_3}{\partial B_n} & \hdots & \frac{\partial x_n}{\partial B_n} & \frac{\partial y_n}{\partial B_n} & \frac{\partial x_2}{\partial B_n} & \frac{\partial y_2}{\partial B_n} \smallskip \\
  \frac{\partial x_1}{\partial B_1} & \frac{\partial y_1}{\partial B_1} & \frac{\partial x_3}{\partial B_1} & \frac{\partial y_3}{\partial B_1} & \hdots & \frac{\partial x_n}{\partial B_1} & \frac{\partial y_n}{\partial B_1} & \frac{\partial x_2}{\partial B_1} & \frac{\partial y_2}{\partial B_1} \smallskip \\
  \frac{\partial x_1}{\partial B_2} & \frac{\partial y_1}{\partial B_2} & \frac{\partial x_3}{\partial B_2} & \frac{\partial y_3}{\partial B_2} & \hdots & \frac{\partial x_n}{\partial B_2} & \frac{\partial y_n}{\partial B_2} & \frac{\partial x_2}{\partial B_2} & \frac{\partial y_2}{\partial B_2}
  \end{pmatrix}\text{.}
\end{align*}\smallskip
We write it as a block matrix with $2\times 2$ blocks:
\begin{align*}
 \tilde{\tilde{M}}=\begin{pmatrix}
  C_1 & C_3 & C_4 &\hdots  & C_n & C_2 \\
  \text{\large 0} & P_3 & \text{\large 0} & \hdots & \text{\large 0} & L_3\\
  \text{\large 0} & \text{\large 0} & P_4 & \hdots & \text{\large 0} & L_4\\
  \vdots & \vdots & \vdots & \ddots & \vdots & \vdots\\
  \text{\large 0} & \text{\large 0} & \hdots & \hdots & P_n & L_n \\
  K_1 & K_3 & \hdots & \hdots & K_n & L_1
   \end{pmatrix}\text{,}
\end{align*}
where $C_1$ is given by \eqref{MatrC1}, the matrices $C_j$ are given by \eqref{MatrCj}, with $A_2$ as in \eqref{A2subs}, and the matrices $P_k$ are given by \eqref{MatrPj}. Moreover,
\begin{align*}
 K_1 = \begin{pmatrix}
        \frac{\partial x_1}{\partial B_1} & \frac{\partial y_1}{\partial B_1} \smallskip \\
        \frac{\partial x_1}{\partial B_2} & \frac{\partial y_1}{\partial B_2}
       \end{pmatrix} =
       \begin{pmatrix}
        -\tfrac{\cosh(tr)-1}{r} & \tfrac{\sinh(rt)}{r} \smallskip\\
        0 & 0
       \end{pmatrix}\text{,}
\end{align*}
\begin{align*}
 K_3 = \dots = K_n = \begin{pmatrix}
                     0 & 0 \\
                     0 & 0
                    \end{pmatrix}\text{,}
\end{align*}

\begin{align*}
 L_1 = \begin{pmatrix}
        \frac{\partial x_2}{\partial B_1} & \frac{\partial y_2}{\partial B_1} \smallskip \\
        \frac{\partial x_2}{\partial B_2} & \frac{\partial y_2}{\partial B_2}
       \end{pmatrix}=
       \begin{pmatrix}
        -\tfrac{B_1\sin(rt)}{A_2r} & -\tfrac{B_1(\cos(rt) - 1)}{A_2r} \smallskip \\
        -\tfrac{B_2\sin(rt)+A_2(\cos(rt)-1)}{A_2r} & \tfrac{A_2\sin(rt) - B_2(\cos(rt)-1)}{A_2r}
       \end{pmatrix}\text{.}
\end{align*}
Finally, similarly as above, we write $\tilde{\tilde{M}}$ in a block form
\begin{align}
 \tilde{\tilde{M}} = \begin{pmatrix}
      P & R \\
      K & L_1
     \end{pmatrix}\text{,}
\end{align}
where this time
\begin{align*}
P=\begin{pmatrix}
  C_1 & C_3 & \hdots & C_n\\
  \text{\large 0} & P_3 & \hdots & \text{\large 0}\\
  \vdots & \vdots & \ddots & \vdots & \\
  \text{\large 0} & \text{\large 0} & \hdots & P_n
   \end{pmatrix}
\end{align*}
is of order $2n-2\times 2n-2$,
\begin{align*}
 R = \begin{pmatrix}
      C_2 \\
      L_3 \\
      \vdots \\
      L_n
     \end{pmatrix}
\end{align*}
is of order $2n-2\times 2$,
\begin{align*}
 K = \begin{pmatrix}
      K_1 & K_3 & \hdots & K_n
     \end{pmatrix}
\end{align*}
is of order $2\times 2n-2$, and $L_1$ is of order $2\times 2$. Notice that the matrix $P$ is invertible. Indeed, as we already know the matrices $P_k$ are nonsingular, and moreover
\[
 \det C_1 = -\tfrac{(\cosh(rt)-1)(B_1^2-1)}{r^2}\ne 0 \text{.}
\]
In fact, its inverse is
\begin{align*}
P^{-1}=\begin{pmatrix}
  C_1^{-1} & -C_1^{-1}C_3P_3^{-1} & \text{\large 0} & -C_1^{-1}C_nP_n^{-1}\\
  \text{\large 0} & P_3^{-1} & \text{\large 0} & \text{\large 0}\\
  \vdots & \vdots & \ddots & \vdots & \\
  \text{\large 0} & \text{\large 0} & \hdots & P_n^{-1}
   \end{pmatrix} \text{.}
\end{align*}
As a result, we can use Schur's formula to compute the determinant of $\tilde{\tilde{M}}$ which gives
\begin{gather*}
 \det\tilde{\tilde{M}} = \det P\det\big(L_1 - KP^{-1}R\big) = \\
 \det C_1(\det P_2)^{n-2}\det\big(L_1 - K_1C_1^{-1}(C_2 - \sum_kC_kP_k^{-1}L_k)\big)\text{.}
\end{gather*}
Proceeding as above, we obtain the formula for $\det\tilde{\tilde{M}}$:
\begin{align*}
 \frac{(-2)^{n-1}(B_1^2-1)\big(\sinh(rt)(\cos(rt)-1)+\sin(rt)(\cosh(rt)-1)\big)(\cos(rt)-1)^{n-2}}{r^{2n}\sqrt{1-B_1^2-B_2^2-\sum_k(A_k^2+B_k^2)}} \text{.}
\end{align*}
This demonstrates that the Jacobian of $\tilde{\tilde{F}}$ is nonzero on $\tilde{\tilde{\mathcal{D}}}$. 

\smallskip

The computations above together with Proposition \ref{propositionProj} show that $W$ is an analytic hypersurface. 

We now turn to the following important result, which is also known in the classical case.
\begin{proposition}\label{NonMax}
 Let $\gamma_i\colon[0,T_i]\longrightarrow\R^{2n+1}$, $i=1,2$, be nspc.f.d. geodesics with the same initial point $\gamma_1(0)=\gamma_2(0)$. Suppose, moreover, that
 \begin{enumerate}
  \item $\gamma_1(t_1) = \gamma_2(t_2)$ for some $t_1\in(0,T_1)$ and $t_2\in(0,T_2)$;
  \item the restriction $\gamma_{i|[0,t_i]}$ is a maximizing geodesic for $i=1,2$;
  \item $\dot{\gamma}_1(t_1)$ and $\dot{\gamma}_2(t_2)$ are not colinear.
 \end{enumerate}
Then, for each $i=1,2$ and every $s_i\in(t_i,T_i]$, the restriction $\gamma_{i|[0,s_i]}$ is not a maximizing geodesic.
\end{proposition}
\begin{proof}
 Suppose to the contrary, that for instance $\gamma_{1|[0,s_1]}$ is a maximizing geodesic. Then the concatenation of $\gamma_{2|[0,t_2]}$ and $\gamma_{1|[t_1,s_1]}$ would be a nonsmooth maximizer. This is impossible, however, since all maximizers are smooth in our setting.
\end{proof}
Let us remark that the proposition does not hold for contact sub-Lorentzian structures on $\R^3$; see \cite{GroJMSc}. See also \cite{Gro2011,Gro2012} for the noncontact case. 
\smallskip

Now we are in a position to describe explicitly the structure of the boundary $\partial J^+(0)$ of the reachable set. Denote the image under the map in \eqref{Map} by $W$. 
Suppose that 
\[
 \gamma\colon[0,+\infty)\longrightarrow\R^{2n+1}
\]
is a null f.d. geodesic starting from $\gamma(0)=0$, with the initial covector as in \eqref{TheSetD}. It follows from the general theory (see \cite{Gro2009}) that the restriction $\gamma_{|[0,T]}$ is geometrically optimal, or equivalently, a null maximizing geodesic, provided that $T>0$ is sufficiently small. On the other hand, the computations above give
\[
 \gamma((0,\tfrac{2\pi}{|r|}))\subset W\text{.}
\]
Consequently, Proposition \ref{PropBound} implies that 
\[
 \gamma_{|(0,\frac{2\pi}{|r|})}
\]
is a trajectory of $\nabla_H\Phi^{\mathrm{sgn\,}r}$, where $\Phi^{\mathrm{sgn\,}r}$ is analytic and $\nabla_H\Phi^{\mathrm{sgn\,}r}$ is a null f.d. vector field. This, in turn, implies that $\gamma$ is geometrically optimal when restricted to $[0,\tfrac{2\pi}{|r|}]$. Furthermore, Proposition \ref{NonMax}, together with the considerations at the beginning of this subsection, shows that $\gamma_{|[0,T]}$ is no longer geometrically optimal for $T>\frac{2\pi}{|r|}$. 
Indeed, if $\gamma$, $\tilde{\gamma}$ have the initial covectors \eqref{InCov1} and \eqref{InCov2}, respectively, where 
\[
 \sum\nolimits_j(A_j^2+B_j^2)=\sum\nolimits_j(\tilde{A}_j^2+\tilde{B}_j^2)
\]
and
\[
 (A_2,B_2,\dots,A_n,B_n) \ne (\tilde{A}_2,\tilde{B}_2,\dots,\tilde{A}_n,\tilde{B}_n) \text{,}
\]
then
\[
 \dot{\gamma}(\tfrac{2\pi}{|r|}) \ne \dot{\tilde{\gamma}}(\tfrac{2\pi}{|r|})\text{.}
\]
Thus, once $t>\frac{2\pi}{|r|}$, the geodesic $\gamma$ remains in $I^+(0)$.

\section{Conjugate and cut locus of the origin}
We will compute the Jacobian of the exponential map with pole at the origin along nspc.f.d. Hamiltonian geodesics; cf. Subsection \ref{subsgen}. To this end, we may take the initial covectors to be of the form 
\[
 (-1,B_1,A_2,B_2,\dots,A_n,B_n,r)\text{.}
\]
That is, we impose the normalization $A_1=-1$ and, consequently, allow $t$ to vary, Thus, at the beginning,  we consider the map
\begin{align*}
 E\colon \mathcal{D}\ni(B_1,A_2,B_2,\dots,A_n,B_n,r,t)\longrightarrow (x_1(t),y_1(t),\dots,x_n(t),y_n(t),z(t))\text{,}
\end{align*}
where
\[
 \mathcal{D} = \{(B_1,A_2,B_2,\dots,A_n,B_n,r,t):\; r\ne 0, \, t>0\} \text{.}
\]
As in the preceding section $x_1(t)$, etc., are given by \eqref{geodrow}. The Jacobian of $E$, up to a sign, can be obtained by computing $\det M$, where
\begin{align*}
 M = \begin{pmatrix}
 \frac{\partial x_1}{\partial r} & \frac{\partial y_1}{\partial r} & \frac{\partial x_2}{\partial r} & \frac{\partial y_2}{\partial r} & \hdots & \frac{\partial x_n}{\partial r} & \frac{\partial y_n}{\partial r} & \frac{\partial z}{\partial r} \smallskip \\ 
  \frac{\partial x_1}{\partial t} & \frac{\partial y_1}{\partial t} & \frac{\partial x_2}{\partial t} & \frac{\partial y_2}{\partial t} & \hdots & \frac{\partial x_n}{\partial t} & \frac{\partial y_n}{\partial t} & \frac{\partial z}{\partial t} \smallskip \\
  \frac{\partial x_1}{\partial A_2} & \frac{\partial y_1}{\partial A_2} & \frac{\partial x_2}{\partial A_2} & \frac{\partial y_2}{\partial A_2} & \hdots & \frac{\partial x_n}{\partial A_2} & \frac{\partial y_n}{\partial A_2} & \frac{\partial z}{\partial A_2} \smallskip \\
  \frac{\partial x_1}{\partial B_2} & \frac{\partial y_1}{\partial B_2} & \frac{\partial x_2}{\partial B_2} & \frac{\partial y_2}{\partial B_2} & \hdots & \frac{\partial x_n}{\partial B_2} & \frac{\partial y_n}{\partial B_2} & \frac{\partial z}{\partial B_2} \smallskip \\
  \vdots & \vdots & \vdots & \vdots & \ddots & \vdots & \vdots & \vdots \smallskip \\
  \frac{\partial x_1}{\partial A_n} & \frac{\partial y_1}{\partial A_n} & \frac{\partial x_2}{\partial A_n} & \frac{\partial y_2}{\partial A_n} & \hdots & \frac{\partial x_n}{\partial A_n} & \frac{\partial y_n}{\partial A_n} & \frac{\partial z}{\partial A_n}\smallskip \\
  \frac{\partial x_1}{\partial B_n} & \frac{\partial y_1}{\partial B_n} & \frac{\partial x_2}{\partial B_n} & \frac{\partial y_2}{\partial B_n} & \hdots & \frac{\partial x_n}{\partial B_n} & \frac{\partial y_n}{\partial B_n} & \frac{\partial z}{\partial B_n} \smallskip \\
  \frac{\partial x_1}{\partial B_1} & \frac{\partial y_1}{\partial B_1} & \frac{\partial x_2}{\partial B_1} & \frac{\partial y_2}{\partial B_1} & \hdots & \frac{\partial x_n}{\partial B_1} & \frac{\partial y_n}{\partial B_1} & \frac{\partial z}{\partial B_1} \smallskip \\
  \end{pmatrix}\text{.}
\end{align*}\smallskip
We write $M$ as a block matrix:
\begin{align*}
 M=\begin{pmatrix}
  C_1 & C_2  &\hdots  & C_n & R_1 \\
  \text{\large 0} & P_2 & \text{\large 0} & \hdots & R_2\\
   \vdots & \vdots & \ddots & \vdots & \vdots\\
  \text{\large 0} & \text{\large 0} & \hdots & P_n & R_n \\
  T_1 & T_2 & \hdots & T_n & S
   \end{pmatrix}\text{,}
\end{align*}
where $C_1$ is given by \eqref{MatrC1}, the matrices $C_j$ are given by \eqref{MatrCj}, and the matrices $P_j$ are given by \eqref{MatrPj}. Moreover,
\begin{gather*}
 R_1 = \begin{pmatrix}
        \frac{\partial z}{\partial r}\smallskip \\
        \frac{\partial z}{\partial t}
       \end{pmatrix}=\\
       \begin{pmatrix}
        \frac{-rt\sum_j(A_j^2+B_j^2)\cos(rt) + 2\sum_j(A_j^2+B_j^2)\sin(rt) + (1-B_1^2)rt\cosh(rt) - 2(1 - B_1^2)\sinh(rt) - (B_1^2 + \sum_j(A_j^2+B_j^2) - 1)rt}{2r^3} \smallskip \\
        \tfrac{-\sum_j(A_j^2+B_j^2)\cos(rt) + (1-B_1^2)\cosh(rt) + B_1^2 + \sum_j(A_j^2+B_j^2) - 1}{2r}
       \end{pmatrix}\text{,}
\end{gather*}
\begin{gather*}
 R_j = \begin{pmatrix}
        \frac{\partial z}{\partial A_j}\smallskip \\
        \frac{\partial z}{\partial B_j}
       \end{pmatrix}=
       \begin{pmatrix}
        \frac{A_j(rt - \sin(rt))}{r^2} \smallskip \\
        \frac{B_j(rt - \sin(rt))}{r^2}
       \end{pmatrix}\text{,}
\end{gather*}
\begin{gather*}
 S = \begin{pmatrix}
        \frac{\partial z}{\partial B_1}
       \end{pmatrix}=
       \begin{pmatrix}
        \frac{B_1(rt - \sinh(rt))}{r^2}
       \end{pmatrix}\text{,}
\end{gather*}
and, finally,
\begin{gather*}
 T_1 = \begin{pmatrix}
        \frac{\partial x_1}{\partial B_1} & \frac{\partial y_1}{\partial B_1}
       \end{pmatrix}=
       \begin{pmatrix}
        -\frac{\cosh(rt)-1}{r} & \frac{\sinh(rt)}{r}
       \end{pmatrix}\text{,}
\end{gather*}
\begin{gather*}
 T_2 = \dots = T_n = \begin{pmatrix}
                      0 & 0
                     \end{pmatrix}\text{.}
\end{gather*}

We now rewrite $M$ as
\begin{align}
 M = \begin{pmatrix}
      P & R \\
      T & S
     \end{pmatrix}\text{,}
\end{align}
where 
\begin{align*}
P=\begin{pmatrix}
  C_1 & C_2 & \hdots & C_n\\
  \text{\large 0} & P_2 & \hdots & \text{\large 0}\\
  \vdots & \vdots & \ddots & \vdots & \\
  \text{\large 0} & \text{\large 0} & \hdots & P_n
   \end{pmatrix}
\end{align*}
is of order $2n\times 2n$ and by \eqref{DetPj} is invertible for $t\ne\frac{2k\pi}{|r|}$, $k=1,2,\dots$,
\begin{align*}
 R = \begin{pmatrix}
      R_1 \\
      \vdots \\
      R_n
     \end{pmatrix}
\end{align*}
is of order $2n\times 1$,
\begin{align*}
 T = \begin{pmatrix}
      T_1 & T_2 & \hdots & T_n
     \end{pmatrix} = \begin{pmatrix}
      T_1 & 0 & \hdots & 0
     \end{pmatrix}
\end{align*}
is of order $1\times 2n$, and $S$ is of order $1\times 1$. 
Suppose $t\ne\frac{2k\pi}{|r|}$. Using Schur's formula we obtain
\begin{gather*}
 \det M = \det P \det (S-TP^{-1}R) = \\
 \det C_1 \cdot(\det P_2)^{n-1} \cdot\big(S - T_1C_1^{-1}(R_1 - \sum\nolimits_jC_jP_j^{-1}R_j)\big) \text{.}
\end{gather*}
After computations we are led to
\begin{gather*}
 \det M = \frac{(-2)^{n-1}(\cos(rt)-1)^{n-2}W(rt)}{r^{2n+2}}
\end{gather*}
with
\begin{gather}
 W(rt) =
 (C-D)rt(\cos(rt) - 1)\sinh(rt) + (\cosh(rt)-1)(2D(\cos(rt) - 1) + Crt\sin(rt)) \text{,}
\end{gather}
where for simplicity we set $C=\sum\nolimits_j(A_j^2+B_j^2)$ and $D = -1 + B_1^2 + C$. Note that the formula for $\det M$ remains valid at
\[
 t_k = \frac{2k\pi}{|r|} \text{,}
\]
$k=1,2,\dots$ Indeed, infinitely many geodesics starting from the origin meet at these values of $t$, as can be seen by an argument similar to that used at the beginning of Subsection \ref{CompCharGeomOpt}. Therefore $\det M_{|t=t_k}=0\frac{}{}$.

Let us notice that $W$ is an even function of $x=rt$, so in order to determine the zeros of $W$, it is enough to consider the expression $W(x)$ for $x>0$. Let us write
\[
 W(x) = C\Pi_1(x) + D\Pi_2(x)\text{,}
\]
where 
\[
 \Pi_1(x) = x\big(\sinh x(\cos x-1) + \sin x(\cosh x -1)\big)\text{,}
\]
\[
 \Pi_2(x) = \big(2(\cosh x - 1) - x\sinh x\big)(\cos x - 1) \text{.}
\]
We now examine the zeros of the two expressions above for $x>0$. Obviously, $\Pi_1(2k\pi)=0$ for every positive integer $k$. On the other hand, for $x\ne 2k\pi$, $\Pi_1$ can be rewritten as
\[
 \Pi_1(x) = x(1-\cos x)(\cosh x-1)(\cot \tfrac{x}{2} - \coth \tfrac{x}{2})\text{.}
\]
We check that the function
\[
 g(x) = \cot \tfrac{x}{2} - \coth \tfrac{x}{2},
\]
is strictly decreasing on every interval $(2k\pi,2(k+1)\pi)$. Since $\lim_{x\to 2k\pi^+}g(x)=+\infty$ and $g(2k\pi+\frac{\pi}{2})<0$, we find that for every $k=1,2,\dots$ there exits a unique number $\alpha_k\in(2k\pi, 2k\pi+\frac{\pi}{2})$ which is also a zero for $\Pi_1$. Thus, $\Pi_1(x)=0$ if and only if $x=2k\pi$ or $x=\alpha_k$, $k=1,2,\dots$.
Next, writing $\Pi_2$ in the form
\[
 \Pi_2(x) = (\cos x-1)(2\tanh\tfrac{x}{2}-x)\sinh x
\]
we see that $\Pi_2(x)\geq 0$, with equality if and only if $x=2k\pi$, $k=1,2\dots$

Finally, we examine the expression $C\Pi_1(x) + D\Pi_2(x)$ under the assumption $D<0$ and $C>0$. For $x\ne 2k\pi$, it can be written as
\[
 x(1-\cos x)(\cosh x-1)\Big[C(\cot \tfrac{x}{2} - \coth \tfrac{x}{2}) + D(\coth\tfrac{x}{2} - \tfrac{2}{x})\Big] \text{.}
\]
Since the first three factors are strictly positive, it is enough to investigate
\[
 h(x) = C(\cot \tfrac{x}{2} - \coth \tfrac{x}{2}) + D(\coth\tfrac{x}{2} - \tfrac{2}{x}) \text{.}
\]
A direct computation shows that $h'(x)<0$ on every interval $(2k\pi,2(k+1)\pi)$, that $\lim_{x\to 2k\pi^+}h(x)=+\infty$ and that $W(\alpha_k)<0$. It follows that for every positive integer $k$ there exists a unique $\beta_k\in(2k\pi,\alpha_k)$ such that $W(\beta_k)=0$. Thus, when $D<0$ and $C>0$, $W(x)=0$ if and only if $x=2k\pi$ or $x=\beta_k$, $k=1,2,\dots$.

Now we consider various cases according to the signs of $C$ and $D$:
\begin{enumerate}
 \item If $D<0$ and $C>0$, then $W(rt)=0$ for $t=\frac{2k\pi}{|r|}$ or $t=\frac{\beta_k}{|r|}$, $k=1,2,\dots$ This case corresponds to t.f.d. geodesics distinct from $t\longrightarrow (t,B_1t,0,\dots,0)$, $|B_1|<1$;
 \item If $D<0$ and $C=0$, then $W(rt)=0$ for $t=\frac{2k\pi}{|r|}$, $k=1,2,\dots$. This case corresponds to t.f.d. geodesics of the form $t\longrightarrow (t,B_1t,0,\dots,0)$, $|B_1|<1$;
 \item If $D=0$ and $C>0$, then $W(rt)=0$ for $t=\frac{2k\pi}{|r|}$ or $t=\frac{\alpha_k}{|r|}$, $k=1,2,\dots$. This case corresponds to null f.d. geodesics distinct from $t\longrightarrow (t,\pm t,0,\dots,0)$;
 \item If $D=0$ and $C=0$, then $W(rt)=0$ for every $t>0$. This case corresponds to the two null f.d. geodesics $t\longrightarrow (t,\pm t,0,\dots,0)$.
\end{enumerate}

As was emphasized above, the Jacobian is smooth at $r=0$ as well, and we may verify that
\[
 \det M_{|r=0} = \tfrac{1}{12}(2C-D)t^{2n+2}\text{.}
\]
For nspc.f.d. geodesics $C\geq 0$ and $D\leq 0$; hence, this expression vanishes if and only if $C=D=0$.
\bigskip

The arguments presented above allow us to write in the explicit form the future null and timelike conjugate loci of the origin, as it is stated in Theorems \ref{Th4} and \ref{Th5}. 
The description of the null cut locus follows from Subsection \ref{CompCharGeomOpt}.

Let us now examine timelike geodesics. First, arguing as at the beginning of Subsection \ref{subssmooth}, we conclude that the t.f.d. geodesics starting from the origin and given by \eqref{GeodEqsrzero} are globally length maximizing.

Next, suppose that 
\[
 \gamma\colon[0,+\infty)\longrightarrow\R^{2n+1}
\]
is a t.f.d. geodesic given by \eqref{GeodEqsrnotzero}, with $A_1=-1$ and $|B_1|<1$. It is known (\cite{Gro2002}) that $\gamma_{|[0,T]}$ is length maximizing for all sufficiently small $T>0$. The above computation of $\det M$ shows that the restriction
\[
 \gamma_{|[0,\frac{2\pi}{|r|})}
\]
is the unique maximizer joining its endpoint. On the other hand, Proposition \ref{NonMax} implies that $\gamma_{|[0,T]}$ is no longer maximizing for $T>\frac{2\pi}{|r|}$ . 

This completes the proof of Theorems \ref{Th4} and \ref{Th5}.

\end{document}